\documentclass[12pt]{amsart}
\usepackage[a4paper,margin=27mm]{geometry}
\usepackage[T1]{fontenc}
\usepackage{xcolor}
\usepackage{lmodern,mathtools,amssymb,mathrsfs,microtype}
\usepackage[colorlinks=true,linkcolor=blue,citecolor=blue,urlcolor=blue]{hyperref}
\hypersetup{pdftitle={A phase transition and critical persistence for superprocesses in random environments}}
\allowdisplaybreaks
\numberwithin{equation}{section}
\newtheorem{theorem}{Theorem}[section]
\newtheorem{lemma}[theorem]{Lemma}
\newtheorem{proposition}[theorem]{Proposition}

\theoremstyle{definition}

\theoremstyle{remark}

\newcommand{\E}{\mathbb E}
\newcommand{\Pp}{\mathbb P}
\newcommand{\R}{\mathbb R}
\newcommand{\Z}{\mathbb Z}

\newcommand{\cH}{\mathcal H_q}
\newcommand{\Mloc}{\mathcal M_{\mathrm{loc}}(\R^d)}
\newcommand{\Mf}{\mathcal M_{\mathrm{f}}(\R^d)}
\newcommand{\Mp}{\mathcal M_{\mathrm{p}}(\R^d)}
\newcommand{\Mvarpi}{\mathcal M_{\varpi}(\R^d)}
\newcommand{\ip}[2]{\langle #1,#2\rangle}

\title[Phase transition and critical persistence]{A phase transition and critical persistence\protect\\ for
 superprocesses in random environments}
\date{}

\author{Jieliang Hong}
\address{Department of Mathematics, Southern University of Science and Technology,
	Shenzhen 518055, China}
\email{hongjl@sustech.edu.cn}

\author{Jie Xiong}
\thanks{This work was supported by the National Key R\&D Program of China (2022YFA1006102),
and the National Natural Science Foundation of China (12471418, 12595294, 12231002, 12571150), and the New Cornerstone Science Foundation (NCI202501).}
\address{Department of Mathematics and Shenzhen International
	Center for Mathematics, Southern University of Science and Technology,
	Shenzhen 518055, China}
\email{xiongj@sustech.edu.cn}

\author{Wen Xu}
\address{School of Mathematical Sciences, Peking University,
	Beijing 100871, China}
\email{xuwen@math.pku.edu.cn}
\usepackage{etoolbox}
\makeatletter
\renewenvironment{proof}[1][\proofname]{%
   \par\pushQED{\qed}\normalfont%
   \topsep6\p@\@plus6\p@\relax
   \trivlist\item[\hskip\labelsep\bfseries#1\@addpunct{.}]%
   \ignorespaces
}{%
   \popQED\endtrivlist\@endpefalse
}
\makeatother

\begin{document}

	
	
\begin{abstract}
We consider a critical superprocess on $\R^d$, $d\ge3$, in a Gaussian environment that is white in time and has spatial covariance $\kappa(1+|x-y|^2)^{-\alpha/2}$, with $\alpha>2$. Starting from Lebesgue measure, the process persists for $0\le\kappa\le\kappa_c$, where $\kappa_c\in(0,\infty)$. For $\kappa>\kappa_c$, each fixed compact set almost surely carries no mass from some finite time onwards. The last time of positive mass in such a set has an exponential tail bound whose decay rate can be chosen independently of the set. By linearizing the limiting Laplace functional at zero, we show that the limiting law of the superprocess has mean measure $\E[Z_\infty^\kappa] dx$, where $Z^\kappa$ is the associated Brownian polymer martingale. For this polymer, the thresholds for strong disorder and strictly negative normalized free energy both equal $\kappa_c$. Weak disorder holds at $\kappa_c$, which lies strictly above the second-moment threshold.
\end{abstract}

\maketitle

\section{Introduction}\label{sec:introduction}

Superprocesses with random branching environments were constructed by
Mytnik~\cite{Myt96} as scaling limits of branching particle systems and
characterized by a martingale problem. We are primarily concerned with the case in which the initial measure is Lebesgue measure, denoted by $m$. Since this initial
measure has infinite total mass, the question concerns the fate of the
population in fixed bounded regions. We use two precise notions of
long-time behaviour. Persistence means that the law of $X_t$, averaged
over both branching and environmental randomness, converges weakly on the
space of locally finite measures equipped with the vague topology to a
nontrivial limiting law whose mean measure is still $m$.
Local extinction means that, for every fixed compact set $K$,
there is an almost surely finite random time $T_K$ such that $X_t(K)=0$
for all $t\ge T_K$. The time $T_K$ may depend on $K$. 

Mytnik and Xiong~\cite{MX} developed
the conditional log-Laplace description and used it to study local
extinction. For bounded covariances satisfying a positive lower bound of order
$1\wedge|x-y|^{-\gamma}$, they proved  local extinction for
$0\le\gamma<2$ and convergence to zero in probability on compact sets for
$\gamma=2$. These conclusions hold in every dimension. For faster spatial
decay in dimensions at least three, their Conjecture~1.4 predicted
persistence at sufficiently small values of the environmental amplitude $\kappa$.  Conditional log-Laplace 
technique was first derived by Xiong \cite{Xio} for another superprocess model with random environment
which affects the spatial motion of the particles. That model was first studied by Skoulakis and Adler \cite{SA}
while leaving the establishment of the conditional log-Laplace equation as an open problem.

Chen, Ren and Zhao~\cite{CRZ} proved the small-environment conjecture under
a quantitative Green-potential condition and established  local
extinction at sufficiently large amplitudes for a fixed smooth covariance
profile. More precisely, their Theorems~1.3 and~1.5 provide the two
one-sided conclusions, while Remark~1.6(ii) asks whether an
extinction--persistence dichotomy holds. Fix \(d\ge 3\) and \(\alpha>2\), and consider the family of spatial covariances
\[
g_\kappa(x,y)=\kappa q(x-y),\qquad q(z)=(1+|z|^2)^{-\alpha/2}.
\]
For the process started from Lebesgue measure, let
\[
\kappa_{\mathrm p}:=\sup\{\kappa\ge0:\text{persistence holds}\},
\qquad
\kappa_{\mathrm e}:=\inf\{\kappa\ge0:\text{local extinction holds}\}.
\]
The two sufficient conditions leave open whether \(\kappa_{\mathrm p}=\kappa_{\mathrm e}\) and, if so, which behaviour occurs at the common threshold.

For this covariance family, we prove that $\kappa_{\mathrm p}=\kappa_{\mathrm e}=\kappa_c$ for some $\kappa_c\in(0,\infty)$ and that the critical point belongs to the persistence phase. Thus, for every $0\le\kappa\le\kappa_c$, the limiting law has mean measure $m$. For every $\kappa>\kappa_c$, each fixed compact set almost surely carries no mass from some finite time onwards. For related results on occupation measures, see~\cite{CHY} for densities in low dimensions and~\cite{FHX} for time averages when $X_0=m$.

The linear equation in the same environment relates the superprocess to a normalized Brownian polymer martingale $Z_t^\kappa$. We write $b_*$ for the threshold for bounded second moments of this martingale. Dawson and Salehi~\cite{DS80} studied persistence for related spatially homogeneous linear evolutions.

Rovira and Tindel~\cite{RT} introduced this Brownian polymer model and gave a second-moment criterion for weak disorder. Lacoin~\cite[Theorem~1.2(i)]{Lacoin11} subsequently considered spatial covariances satisfying $Q(x)\asymp |x|^{-\alpha}$ at infinity. For $d\ge3$ and $\alpha>2$, he proved $0<\beta_c\le\bar\beta_c$, where $\beta_c$ is the weak/strong-disorder threshold in inverse temperature and $\bar\beta_c$ is the threshold for strictly negative normalized free energy. With $Q=q$ his inverse temperature is $\beta=\sqrt\kappa$. For this profile we obtain
\[
 \sqrt{b_*}<\beta_c=\bar\beta_c=\sqrt{\kappa_c},
\]
and weak disorder holds at $\beta_c$. The decay assumption matters for the equality of the latter two thresholds. For radial correlations of order $|x|^{-2}$ in $d\ge3$, Cao, Rang and Wu~\cite[Theorem~1.2(i)]{CRW} showed that strong disorder holds at every positive inverse temperature while the free energy remains zero at sufficiently small inverse temperatures.

For a Brownian polymer defined using a compactly supported spatial mollifier, Bazaes, Lammers and Mukherjee~\cite[footnote~1]{BLM} state that the weak-disorder threshold is strictly larger than the second-moment threshold. Here the strict separation holds for the polynomially decaying covariance $q$, including the range $2<\alpha\le d$ in which $q$ is not integrable over $\R^d$.

The link to the superprocess comes from an identity for the mean of its limiting law. Chen, Ren and Zhao~\cite[Lemma~4.4]{CRZ} established existence of the limit from Lebesgue measure and showed that its mean measure is bounded above by $m$. We identify the mean measure for every $\kappa$ as $m$ multiplied by the expected terminal value of the polymer martingale. To prove this identity we linearize the nonlinear Laplace functional at zero after taking the long-time limit. A pathwise space--time energy estimate justifies this step even when the polymer has unbounded second moments.

For lattice directed polymers in i.i.d.\ environments, Junk and Lacoin~\cite[Theorem~2.1]{JL} proved coincidence of the strong- and very-strong-disorder thresholds under their hypotheses. Their excursion and size-bias argument guides our proof of exponential fractional-moment decay in continuous space. A positive eigenfunction at $b_*$ gives the auxiliary function for the excursion estimate. A comparison on a unit cube controls the supremum over the polymer's starting point and allows us to deduce eventual local extinction for the superprocess started from $m$.

Birkner~\cite{Birkner} obtained a weak-disorder criterion for lattice polymers by size-biasing the partition function. Here we condition the size-biased partition function on the noise along a Brownian path. To prove the strict inequality $b_*<\kappa_c$, we subtract a small rank-one kernel from the time-one Feynman--Kac kernel at $b_*$. The modified kernel has a bounded positive harmonic function and a local limit with Brownian heat-kernel asymptotics. A rank-one splitting of the type introduced by Nummelin~\cite{Nummelin} yields a renewal law with gap probabilities $K(n)\sim Cn^{-d/2}$. Increasing $\kappa$ from $b_*$ to $b_*+\delta$ changes the weight of each gap by a factor of at most $1+C\delta$, uniformly in its length. Berger and Lacoin's criterion~\cite[Theorem~2.2]{BL} gives a strictly positive critical reward for the corresponding disordered pinning model, including the marginal case $d=3$. This absorbs the extra reward $\log(1+C\delta)$ per contact. Mourrat's summability theorem~\cite{Mourrat,AZ} and Kahane's Gaussian comparison~\cite{Kahane} then give uniform integrability of the polymer martingale for all sufficiently small $\delta>0$.

\section{Preliminaries and main results}\label{sec:main}

We write $m(dx)=dx$ for Lebesgue measure and
$\mu(f)=\ip{f}{\mu}=\int f\,d\mu$ whenever the integral exists.
The space $\Mloc$ consists of nonnegative Radon measures that are finite on
compact sets. It is equipped with the vague topology, generated by the maps
$\mu\mapsto\mu(f)$ for $f\in C_c(\R^d)$. Fix an exponent $r>d$, set
$\varpi(x)=(1+|x|^2)^{-r/2}$, and define
\[
 \Mvarpi=\{\mu\in\Mloc:\mu(\varpi)<\infty\}.
\]
This weighted state space contains $m$ and is used to construct processes with infinite initial mass. When a topology on $\Mvarpi$ is needed,
we use weak convergence of the finite measures $\varpi\mu$ against bounded
continuous test functions. All convergence of probability laws below
is weak convergence on $\Mloc$ with its vague topology. We write $\Mf$ and $\Mp$ for the spaces of finite measures and probability measures on $\R^d$, respectively.

For an open set $O$, $C_c^k(O)$ denotes the space of $k$ times continuously
differentiable functions with compact support in $O$. We use $C_c(O)$ for
continuous functions with compact support, $C_b^k(O)$ for functions whose
derivatives through order $k$ are continuous and bounded, and $C_0(\R^d)$
for continuous functions that vanish at spatial infinity, equipped with
the supremum norm. A superscript $+$ restricts any such class to
nonnegative functions. The spaces $L^p(O)$ and $W^{k,p}(O)$ are defined
using Lebesgue measure and weak derivatives, respectively; $H^k=W^{k,2}$.
The notation $\|\cdot\|_p$ means the $L^p(\R^d)$ norm unless another
measure or domain is displayed. For integrable products,
$\ip{f}{h}=\int fh\,dx$. The rank-one kernel $f\otimes h$ acts by
$(f\otimes h)v(x)=f(x)\int h(y)v(y)\,dy$.
For a bounded operator $A$ on a specified Banach space,
$r(A)=\lim_{n\to\infty}\|A^n\|^{1/n}$ denotes its spectral radius.

Let $\cH$ be the Hilbert-space completion of the linear span of symbols
$e_x$, $x\in\R^d$, with inner product
$\ip{e_x}{e_y}_{\cH}=q(x-y)$, after quotienting out vectors of zero norm.
The space is separable, since $x\mapsto e_x$ is continuous and rational
points are dense in $\R^d$. On a probability space
$(\Omega,\mathcal F,\Pp)$ let $\mathcal W$ be a cylindrical Brownian motion
on $\cH$, and set $W(t,x)=\mathcal W_t(e_x)$. Pairings against
$d\mathcal W_s$ below denote cylindrical stochastic integrals.  Thus
\[
 \E[W(t,x)W(s,y)]=(t\wedge s)q(x-y).
\]
The filtration $(\mathcal F_t^W)$ is the completed, right-continuous natural
filtration of this noise. Stochastic integrals with respect to $W$ are It\^o
integrals. 

We distinguish the two heat kernels
\[
 p_t(x)=(2\pi t)^{-d/2}e^{-|x|^2/(2t)},\qquad
 \rho_t(x)=(4\pi t)^{-d/2}e^{-|x|^2/(4t)}.
\]
They correspond to $\frac12\Delta$ and $\Delta$, respectively. The reference
motion of the superprocess and of a single polymer has generator
$\frac12\Delta$. The difference of two independent reference motions has
generator $\Delta$. Expectations denoted by $E_x^B$ and $E_x^D$ refer to
these two path laws, started at $x$ and independent of the environment.
We write $B(x,R)$ for the open Euclidean ball and $B_R=B(0,R)$.
The symbols $c,C$ denote finite positive constants that may change between
displays. We specify their dependence on parameters when a uniform estimate is needed.

Let $\Pp_\mu^\kappa$ denote the annealed law of the superprocess with
initial measure $\mu\in\Mvarpi$ and covariance $g_\kappa(x,y)=\kappa q(x-y)$, and let $\E_\mu^\kappa$ denote expectation under this law. The superscript $\kappa$ is
omitted when fixed.  We write $\mathcal L_\mu^\kappa(X_t)$ for the distribution of $X_t$
under $\Pp_\mu^\kappa$. The process $X_t$ is
specified by the following martingale problem: for
$f\in C_c^2(\R^d)$,
\[
 M_t(f)=X_t(f)-\mu(f)-\frac12\int_0^tX_s(\Delta f)\,ds
\]
is a continuous local martingale in the joint filtration carrying the
process and its noises, with quadratic variation
\begin{align}\label{eq:martingale-problem}
 [M(f)]_t
 &=\int_0^tX_s(f^2)\,ds+\kappa\int_0^t\iint q(x-y)f(x)f(y)\,X_s(dx)X_s(dy)\,ds.
\end{align}
The process is the solution supplied by the conditional construction in
\cite{MX,CRZ}; in particular, its paths are continuous in the vague
topology.  For the initial condition $X_0=m$, we use the construction with
Lebesgue initial measure given in
\cite[Proposition~4.1 and Appendix~A]{CRZ}. We use the conditional branching property with a common
environment. This means that, for $\mu=\mu_1+\mu_2$, the processes started
from $\mu_1$ and $\mu_2$ are independent conditionally on that environment
and their sum has law $\Pp_\mu^\kappa$. They need not be independent under
the annealed law. For example, if $B$ is an independent Brownian path, then
\[
H_t(B):=\int_0^tW(ds,B_s)
=\int_0^t\ip{e_{B_s}}{d\mathcal W_s}_{\cH}
\]
is well defined and, conditionally on $B$, is a centered Gaussian random variable with variance $t$.

The associated forward equations, on the environmental probability
space, are
\begin{align}\label{eq:dual-equations}
 dv_t(x)&=\tfrac12\Delta v_t(x)\,dt+\sqrt\kappa\,v_t(x)\,dW_t(x), \qquad
       v_0(x)=f(x),\notag\\
 du_t(x)&=(\tfrac12\Delta u_t(x)-\tfrac12u_t^2(x))\,dt
          +\sqrt\kappa\,u_t(x)\,dW_t(x), \qquad
       u_0(x)=f(x).
\end{align}
Here $dW_t(x)$ denotes the It\^o differential in time at the fixed
spatial point $x$. 
For bounded nonnegative initial data, solutions are understood in the mild
sense, using convolution with $p_t$. The spatial stochastic convolution is
defined by the Hilbert-space integral above. For instance, the stochastic
term in the equation for $v_t(x)$ is
\[
 \sqrt\kappa\int_0^t\left\langle
       \int p_{t-s}(x-y)v_s(y)e_y\,dy,\,d\mathcal W_s
                                      \right\rangle_{\cH}.
\]

In the conditional construction, the log-Laplace equation is solved backward
from the observation time. The time-reversed solution has the law of the
forward solution in \eqref{eq:dual-equations}. Averaging over the environment gives
\begin{equation}\label{eq:annealed-laplace}
 \E_\mu^\kappa e^{-X_t(f)}
       =\E e^{-\mu(u_t^f)},\qquad f\in C_c^+(\R^d).
\end{equation}
Section~\ref{sec:bridge} gives the conditional formula and its extension
to bounded test functions for finite $\mu$.

For later reference, we write
\begin{equation}\label{eq:polymer-definition}
 Z_t^\kappa(x)=E_x^B
       \exp\{\sqrt\kappa H_t(B)-\kappa t/2\},\qquad
 Z_t^\kappa=Z_t^\kappa(0).
\end{equation}
It is a nonnegative mean-one martingale in $\mathcal F_t^W$ and therefore
has an almost sure limit $Z_\infty^\kappa$. We call the parameter weakly
disordered when this martingale is uniformly integrable, and strongly
disordered when $Z_\infty^\kappa=0$ almost surely. By the zero--one law for $\{Z_\infty^\kappa=0\}$ in \cite[Section~1.2]{Lacoin11} and Proposition~\ref{prop:mean-dichotomy} below, weak disorder is equivalent to $\Pp(Z_\infty^\kappa>0)=1$.

The first theorem describes the long-time behaviour of the superprocess started from Lebesgue measure.

\begin{theorem}[A common threshold and its endpoint]\label{main}
There exists
$\kappa_c\in(0,\infty)$ such that the following statements hold when
$X_0=m$.
\begin{enumerate}
\item If $0\le\kappa\le\kappa_c$, then
$\mathcal L_m^\kappa(X_t)\Rightarrow\pi_\kappa$ as $t\to\infty$, where
$\pi_\kappa$ is a nontrivial probability measure on $\Mloc$ and
\begin{equation}\label{eq:limiting-mean-intro}
 \int_{\Mloc}\mu(f)\,\pi_\kappa(d\mu)=\int_{\R^d}f(x)\,dx,
 \qquad f\in C_c^+(\R^d).
\end{equation}
\item If $\kappa>\kappa_c$, there exists $c=c(\kappa)>0$ such that the following holds for every compact $K\subset\R^d$. There is an almost surely finite random time $T_K$ such that $X_s(K)=0$ for all $s\ge T_K$ almost surely. Let $L_K:=\sup\bigl(\{s\ge0:X_s(K)>0\}\cup\{0\}\bigr)$ be the last time of positive mass in $K$. The time $T_K$ can be chosen so that, for some $C_{\kappa,K}<\infty$,
\[
 \Pp_m^\kappa(L_K>t)\le \Pp_m^\kappa(T_K>t)\le C_{\kappa,K}e^{-ct},
 \qquad t\ge0.
\]
\end{enumerate}
Moreover,
\[
 \kappa_c=\sup\{\kappa\ge0:(Z_t^\kappa)_{t\ge0}
                                   \text{ is uniformly integrable}\}.
\]
In particular, $\kappa_{\mathrm p}=\kappa_{\mathrm e}=\kappa_c$, and persistence holds at the common threshold.
\end{theorem}

Let $G=(-\Delta)^{-1}$ be the Green operator for the difference motion,
\[
 Gf(x)=\int_0^\infty(\rho_t*f)(x)\,dt
      =\frac1{(d-2)|\mathbb S^{d-1}|}
                  \int_{\R^d}|x-y|^{2-d}f(y)\,dy
\]
whenever the integral is absolutely convergent. Define
\begin{equation}\label{eq:second-moment-threshold-intro}
 Tf=G(qf)\quad\text{on }C_0(\R^d),\qquad
 b_*:=r(T)^{-1}.
\end{equation}
Here $|\mathbb S^{d-1}|$ denotes the surface area of the unit sphere.
The operator $T$ is positive and compact. Proposition~\ref{prop:exact-L2} shows that $\sup_{t\ge0}\E[(Z_t^\kappa)^2]<\infty$ holds exactly when $\kappa<b_*$. The next theorem identifies the behaviour beyond this range.

\begin{theorem}[Weak disorder beyond the second-moment threshold]\label{thm:critical-intro}
The thresholds satisfy
\[
 0<b_*<\kappa_c<\infty.
\]
For every $\kappa\in[b_*,\kappa_c]$, the martingale $(Z_t^\kappa)_{t\ge0}$ is uniformly integrable and
\[
 \E Z_\infty^\kappa=1,\qquad
 \lim_{t\to\infty}\E[(Z_t^\kappa)^2]=\infty.
\]
The free energy $p(\kappa):=\lim_{t\to\infty}t^{-1}\E\log Z_t^\kappa$ vanishes for $0\le\kappa\le\kappa_c$ and is strictly negative for $\kappa>\kappa_c$. At $\kappa=\kappa_c$,
\[
 \lim_{t\to\infty}\frac1t\log Z_t^{\kappa_c}=0
                         \quad\Pp\text{-a.s.},
 \qquad
 \lim_{t\to\infty}\frac1t\log\E[(Z_t^{\kappa_c})^\theta]=0
                         \quad\text{for every }\theta\in(0,1).
\]
\end{theorem}

\medskip
\noindent\textbf{Proof strategy and key identities.}
Section~\ref{sec:bridge} establishes convergence to a limiting law
$\pi_\kappa$ for every $\kappa\ge0$ when the superprocess starts from
Lebesgue measure. For $f\in C_c^{\infty,+}(\R^d)$, let
$V_t^f=\int_{\R^d}v_t^f(x)\,dx$ denote the total mass of the linear
solution, and let $V_\infty^f=\lim_{t\to\infty}V_t^f$ be its almost
sure limit. We prove that the mean measure of $\pi_\kappa$ is
determined by the expected terminal mass of this linear equation.
\begin{equation}\label{eq:intro-contribution-bridge}
	\int_{\Mloc}\mu(f)\,\pi_\kappa(d\mu)
	=\E V_\infty^f
	=c_\kappa\int_{\R^d}f(x)\,dx,
	\qquad c_\kappa\in\{0,1\}.
\end{equation}
Here $c_\kappa=\E Z_\infty^\kappa$ is independent of the test function $f$.
The same coefficient determines uniform integrability for every finite
initial measure of the linear equation. This connects the superprocess
transition to a single polymer parameter. Passing from the finite-time Laplace identity to this conclusion requires taking the small-initial-data and long-time limits together.

The next step is the continuous-space implication
\begin{equation}\label{eq:intro-contribution-sharpness}
 \left.
 \begin{array}{c}
   Z_\infty^\kappa=0\quad\Pp\text{-a.s.}\\
   \kappa>b_*
 \end{array}\right\}
 \quad\Longrightarrow\quad
 \forall\theta\in(0,1),\quad
 \E[(Z_t^\kappa)^\theta]\le C_\theta e^{-c_\theta t}.
\end{equation}
Lemma~\ref{lem:derivatives} gives a fractional-moment bound on the
supremum of the polymer mass over starting points in a unit cube.
We use this bound in the finite-volume criterion of
Lemma~\ref{lem:finitevolume} to prove \eqref{eq:intro-contribution-sharpness},
and in the expanding-ball estimate of Lemma~\ref{lem:expanding-ball}
to deduce local extinction when the process starts from Lebesgue measure.
Together with \eqref{eq:intro-contribution-bridge}, these estimates
give persistence for $\kappa<\kappa_c$ and local extinction for
$\kappa>\kappa_c$.

To determine the behaviour at $\kappa_c$, we first prove that weak disorder extends to the right of the second-moment threshold. More precisely,
\begin{equation}\label{eq:intro-contribution-endpoint}
 \text{there exists }\delta_*>0\text{ such that }\quad
 \E Z_\infty^\kappa=1
 \quad\text{for all }\kappa\in[b_*,b_*+\delta_*].
\end{equation}
Section~\ref{sec:critical} obtains \eqref{eq:intro-contribution-endpoint} by comparing the perturbed gap weights with the critical renewal law. It follows that $\kappa_c\ge b_*+\delta_*>b_*$. Continuity of the free energy gives $p(\kappa_c)=0$. If strong disorder held at $\kappa_c$, \eqref{eq:intro-contribution-sharpness} would give $p(\kappa_c)<0$. Thus the polymer is uniformly integrable at $\kappa_c$, and \eqref{eq:intro-contribution-bridge} gives persistence there.

The rest of this article is organized as follows. Section~\ref{sec:bridge} proves the nonlinear-to-linear mean identity and identifies
the common disorder parameter. Section~\ref{sec:sharpness} proves
fractional-moment sharpness, and Section~\ref{sec:eventual-extinction}
deduces local extinction. Section~\ref{sec:critical} treats the
critical parameter and completes the proofs of the two theorems.

\section{The nonlinear-to-linear bridge}\label{sec:bridge}

Throughout this section, $\kappa\ge0$ is fixed. Expectations
without an initial-measure subscript are taken only over the environment. For
nonnegative initial data $f$, write $u_t^f$ and $v_t^f$ for the nonlinear
and linear solutions of \eqref{eq:dual-equations}, respectively. 
\subsection{Duality, total masses and the limiting law}
Fix $T>0$ and a finite initial measure $\mu$. In the joint construction
of \cite[Theorem~2.15]{MX}, denote by $\psi_{s,T}^{\,g}$ the solution of
the backward log-Laplace equation \cite[(2.4)]{MX} with terminal value
$g$ at time $T$. For every bounded nonnegative Borel function $f$,
\begin{equation}\label{eq:quenched-reversed-laplace}
 \E_\mu\!\left[e^{-X_T(f)}\mid\mathcal F_\infty^W\right]
       =\exp\{-\mu(\psi_{0,T}^{\,f})\}\qquad\text{a.s.}
\end{equation}
Indeed, $f_R=f\mathbf1_{B_R}$, $R\in\mathbb N$, is bounded by a constant
multiple of $p_1$, so \cite[(2.22)]{MX} gives this identity for $f_R$.
Use the same backward noise for every $R$. By comparison the fields
$\psi_{0,T}^{\,f_R}$ increase with $R$. For general $f$, set
$\psi_{0,T}^{\,f}=\lim_{R\to\infty}\psi_{0,T}^{\,f_R}$. Since $1-e^{-X_T(f_R)}\uparrow1-e^{-X_T(f)}$,
conditional monotone convergence gives \eqref{eq:quenched-reversed-laplace}.

Under the time change $r=T-s$, the backward equations for the countable
family $(f_R)_R$ are forward equations driven by one cylindrical Brownian
motion \cite[Remark~2.5 and (2.6)]{MX}. Hence
$(\mu(\psi_{0,T}^{\,f_R}))_R$ has the joint law of
$(\mu(u_T^{f_R}))_R$ with one forward noise. On the forward space,
the difference $u^f-u^{f_R}$ solves the linear equation with initial
value $f-f_R$ and additional killing rate $(u^f+u^{f_R})/2$.
Comparison and the linear first moment give
\[
 0\le\E\,\mu(u_T^f-u_T^{f_R})
 \le\E\,\mu(v_T^{f-f_R})=\mu\bigl(p_T*(f-f_R)\bigr)\longrightarrow0.
\]
Here bounded convergence applies since $\mu$ is finite. Thus
$\mu(u_T^{f_R})\uparrow\mu(u_T^f)$ almost surely, and
$\mu(\psi_{0,T}^{\,f})$ has the law of $\mu(u_T^f)$. Averaging
\eqref{eq:quenched-reversed-laplace} gives \eqref{eq:annealed-laplace}
for every bounded nonnegative Borel $f$ and finite $\mu$.

The annealed Laplace identity \eqref{eq:annealed-laplace}
extends from finite initial measures to
$\mu\in\mathcal M_\varpi$ for $f\in C_c^+$.
This extension follows by approximating $\mu$ with its
restrictions to increasing balls and applying the conditional
branching construction and monotone convergence.
For these initial measures, the first-moment formula reads
\begin{equation}\label{eq:first-moment}
	\E_\mu X_t(f)=\mu(p_t*f),\qquad f\ge0.
\end{equation}
 In particular, Tonelli's theorem gives
\[
 \E_mX_t(f)=\int(p_t*f)(x)\,dx
 =\int f(y)\left(\int p_t(x-y)\,dx\right)dy=m(f).
\]

Let $v^g$ denote the solution of the linear equation in
\eqref{eq:dual-equations} with initial condition $v_0^g=g$.
Let $\mathcal K_{s,t}^\kappa(x,y)$ denote its positive solution kernel, where $x$ is the starting position at time $s$
and $y$ is the terminal position at time $t$.  The kernel admits a jointly measurable version and satisfies,
for each deterministic $0\le s<r<t$, the composition identity
\[
\mathcal K_{s,t}^\kappa(x,y)
=\int_{\R^d}\mathcal K_{s,r}^\kappa(x,z)
\mathcal K_{r,t}^\kappa(z,y)\,dz
\qquad\text{a.s.}.
\]
Then 
\[
v_t^g(y)=\int g(x)\mathcal K_{0,t}^\kappa(x,y)\,dx.
\]
The Feynman--Kac representation gives, for every
Borel set $A\subseteq\R^d$,
\[
\int_A\mathcal K_{0,t}^\kappa(x,y)\,dy
=E_x^B\!\left[
\exp\left\{\sqrt\kappa H_t(B)-\frac{\kappa t}{2}\right\}
\mathbf1_{\{B_t\in A\}}
\right].
\]
Taking $A=\R^d$ yields
\[
\int\mathcal K_{0,t}^\kappa(x,y)\,dy
=E_x^B\exp\left\{\sqrt\kappa H_t(B)-\frac{\kappa t}{2}\right\}
=Z_t^\kappa(x),
\]
where the last equality is precisely \eqref{eq:polymer-definition}.
These identities
first hold for bounded integrable $g$ and extend by positivity to
$g\in L^1_+$. 

Define
\begin{equation}\label{l2}
V_t^g:=\int v_t^g(x)\,dx =\int_{\R^d}g(x)
\left(\int_{\R^d}\mathcal K_{0,t}^\kappa(x,y)\,dy\right)dx=\int_{\R^d}g(x)
Z_t^\kappa(x)dx, \quad g\in L^1_+
\end{equation}
The first two moments of $V_t^g$ satisfy
\begin{equation}\label{eq:linear-mass-moments}
 \E V_t^g=m(g),\qquad
 \E[(V_t^g)^2]\le e^{\kappa t}m(g)^2.
\end{equation}
The flow property and the independence of future noise give, for $0\le s<t$,
\[
 \E[V_t^g\mid\mathcal F_s^W]
 =\int v_s^g(x)\,
       \E\left[\int\mathcal K_{s,t}^\kappa(x,y)\,dy\right]dx
 =\int v_s^g(x)\,dx=V_s^g.
\]
Thus $V^g$ is a continuous nonnegative martingale in
$(\mathcal F_t^W)$ and has a finite almost sure limit $V_\infty^g$.

Fix $f\in C_c^{\infty,+}$ and write $F=m(f)$. For $\lambda>0$, define
\[
 U_t(\lambda)=\int u_t^{\lambda f}(x)\,dx,
 \qquad
 Q(g)=\iint q(x-y)g(x)g(y)\,dx\,dy.
\]
Comparison gives $0\le u_t^{\lambda f}\le\lambda v_t^f$. To integrate
\eqref{eq:dual-equations}, take a smooth radial cutoff
$\chi_R(x)=\chi(x/R)$ increasing to $1$, equal to $1$ on $B_R$
and supported in $B_{2R}$. For $g=v^f$ or $g=u^{\lambda f}$, put
$a=1$ or $a=\lambda$, respectively. Then $0\le g_s\le a v_s^f$, so on
each finite interval $[0,T]$ the expected absolute Laplacian error is at
most $CaTFR^{-2}$. Since $q\le1$, the expected bracket of the discarded
noise is bounded by
\[
 \kappa\int_0^T \E\!\left[\left(\int(1-\chi_R)g_s\,dx\right)^2\right]ds
 \longrightarrow0
\]
by dominated convergence and \eqref{eq:linear-mass-moments}.
For the absorption term, two independent Brownian bridges in the
Feynman--Kac formula contribute a joint exponential moment at most
$e^{\kappa s}$, since their integrated covariance is at most $s$. Thus
$\E[(v_s^f(x))^2]\le e^{\kappa s}(p_s*f)(x)^2$. This is integrable in
$(s,x)$ on $[0,T]\times\R^d$ because $\|p_s*f\|_2\le\|f\|_2$.
Letting $R\to\infty$ in the two weak equations gives
\begin{align}\label{eq:integrated-dual}
 dV_t^f&=\sqrt\kappa\int v_t^f(x)\,W(dt,dx),\notag\\
 dU_t(\lambda)&=-\tfrac12\|u_t^{\lambda f}\|_2^2\,dt
             +\sqrt\kappa\int u_t^{\lambda f}(x)\,W(dt,dx).
\end{align}
In particular, $U(\lambda)$ is a continuous
nonnegative supermartingale. We denote its finite almost sure limit by
$U_\infty(\lambda)$. The first-moment bound \eqref{eq:first-moment} makes the family
$\{\mathcal L_m(X_t):t\ge0\}$ tight in $\mathcal M_{\mathrm{loc}}$.
Indeed, choose $\psi_j\in C_c^+$ with $\psi_j\ge1$ on
$\overline B_j$, and define
\[
 L_j=\varepsilon^{-1}2^j m(\psi_j),\qquad
 \mathcal A_\varepsilon=\{\mu:\mu(\psi_j)\le L_j\text{ for all }j\ge1\}.
\]
This is a closed, uniformly locally bounded set of measures and hence is
vaguely compact. Markov's inequality gives
\[
 \sup_{t\ge0}\Pp_m\{X_t\notin\mathcal A_\varepsilon\}
 \le\sum_{j\ge1}\frac{\E_mX_t(\psi_j)}{L_j}
 =\varepsilon\sum_{j\ge1}2^{-j}=\varepsilon.
\] For each $f\in C_c^{\infty,+}$, the annealed duality and bounded
convergence give
\[
 \lim_{t\to\infty}\E_m e^{-\lambda X_t(f)}
       =\E e^{-U_\infty(\lambda)}.
\]
Laplace functionals on this class determine a probability law on
$\mathcal M_{\mathrm{loc}}$. All subsequential limits therefore agree.
We denote the limit by $\pi_\kappa$. In particular, for every $f\in C_c^{\infty,+}$ and $\lambda>0$,
\begin{equation}\label{eq:pi-laplace}
	\int e^{-\lambda\mu(f)}\,\pi_\kappa(d\mu)
	=\E e^{-U_\infty(\lambda)}.
\end{equation}

\begin{lemma}[Almost sure space--time energy]\label{lem:linear-energy}
For every $f\in C_c^{\infty,+}$, we have
\[
 \int_0^\infty\|v_t^f\|_2^2\,dt<\infty,\qquad\text{a.s.}.
\]
\end{lemma}
\begin{proof}
For every fixed $h>0$, there exists $a>0$ such that
\begin{equation}\label{gaussdom}
 q(z)\ge a p_{2h}(z),\qquad z\in\mathbb R^d.
\end{equation}
Suppose first that $\kappa>0$. From \eqref{eq:integrated-dual},
\[
 [V^f]_t=\kappa\int_0^tQ(v_s^f)\,ds.
\]
Write $A_t=[V^f]_t$. By the Dambis--Dubins--Schwarz theorem, we may enlarge the probability
space if necessary and choose a standard one-dimensional Brownian
motion $(\widehat B_a)_{a\ge0}$ such that
\[
V_t^f=F+\widehat B_{A_t},\qquad t\ge0.
\]  On $\{A_\infty=\infty\}$, continuity
of $A_t$ implies that its range is $[0,\infty)$. Nonnegativity of $V^f$
would therefore imply
\[
 \inf_{a\ge0}\widehat B_a\ge-F.
\]
Almost surely, $\inf_{a\ge0}\widehat B_a=-\infty$,
so $\Pp(A_\infty=\infty)=0$ and
\begin{equation}\label{xu1}
\int_0^\infty Q(v_s^f)\,ds<\infty \qquad\text{a.s.}.
\end{equation}

The future kernel is independent of $\mathcal F_t^W$, and $v_t^f$ is
$\mathcal F_t^W$-measurable. The two-particle Feynman--Kac formula gives
\[
 \E[\mathcal K_{t,t+h}^\kappa(y,x)
       \mathcal K_{t,t+h}^\kappa(z,x)\mid\mathcal F_t^W]
 \le e^{\kappa h}p_h(x-y)p_h(x-z).
\]
Multiplying by $v_t^f(y)v_t^f(z)$ and integrating in $x,y,z$ yields
\begin{align}\label{eq:delayed-energy-bound}
 \E[\|v_{t+h}^f\|_2^2\mid\mathcal F_t^W]
 &\le e^{\kappa h}\iint
             v_t^f(y)v_t^f(z)p_{2h}(y-z)\,dy\,dz\le a^{-1}e^{\kappa h}Q(v_t^f).
\end{align}
For each integer $N\ge1$, let
\[
 \tau_N=\inf\left\{t:\int_0^tQ(v_s^f)ds\ge N\right\},
 \qquad \inf\varnothing=\infty.
\]
Since $\{t<\tau_N\}\in\mathcal F_t^W$, conditional expectation and
Tonelli's theorem give
\begin{align*}
 \E\int_0^\infty\mathbf1_{\{t<\tau_N\}}\|v_{t+h}^f\|_2^2dt
 &=\int_0^\infty\E\left[
   \mathbf1_{\{t<\tau_N\}}
   \E[\|v_{t+h}^f\|_2^2\mid\mathcal F_t^W]\right]dt\\
 &\le a^{-1}e^{\kappa h}\E\int_0^{\tau_N}Q(v_t^f)\,dt\le a^{-1}e^{\kappa h}N.
\end{align*}
Therefore,
all the random variables $\int_0^\infty
\mathbf1_{\{t<\tau_N\}}\|v_{t+h}^f\|_2^2\,dt$ are finite on a common event of
probability one. On each sample path for which the integral in \eqref{xu1} is finite, choose an
integer $N$ larger than its value. Then
\[
\int_0^t Q(v_s^f)\,ds<N
\qquad\forall \ t<\infty,
\]
and hence $\tau_N=\infty$. Consequently,
\[
 \Pp\left(\bigcup_{N\ge1}\{\tau_N=\infty\}\right)
=\Pp\left(\int_0^\infty Q(v_s^f)\,ds<\infty\right)=1.
\]
Thus, on this event of probability one, we obtain
\[
\int_0^\infty\|v_{t+h}^f\|_2^2\,dt
=\int_h^\infty\|v_s^f\|_2^2\,ds< \infty.
\]
On the remaining interval, we have
\[
 \E\int_0^h\|v_t^f\|_2^2dt
 \le\|f\|_2^2\int_0^h e^{\kappa t}dt<\infty.
\]

For $\kappa=0$, Young's inequality gives
\begin{align*}
 \int_0^\infty\|v^f_t\|_2^2dt=  \int_0^\infty\|p_t*f\|_2^2dt
 &\le\|f\|_2^2+\|f\|_1^2\int_1^\infty\|p_t\|_2^2dt\\
 &=\|f\|_2^2+(4\pi)^{-d/2}\|f\|_1^2
                           \int_1^\infty t^{-d/2}dt<\infty,
\end{align*}
where the final integral is finite precisely when $d>2$. This completes the proof.
\end{proof}

\begin{theorem}[The nonlinear-to-linear identity]\label{thm:mean-bridge}
For every $f\in C_c^{\infty,+}$, we have
\[
 \lim_{n\to\infty}nU_\infty(1/n)=V_\infty^f
       \qquad\text{a.s.}.
\]
Moreover,
\begin{equation}\label{meanbridge}
 \int\mu(f)\,\pi_\kappa(d\mu)=\E V_\infty^f.
\end{equation}
\end{theorem}
\begin{proof}
The rescaled solution $w_t^{(n)}:=nu_t^{f/n}$ satisfies
	\[
	dw_t^{(n)}(x)
	=\left(\frac12\Delta w_t^{(n)}(x)
	-\frac1{2n}\bigl(w_t^{(n)}(x)\bigr)^2\right)dt
	+\sqrt\kappa\,w_t^{(n)}(x)\,dW_t(x),
	\qquad w_0^{(n)}=f.
	\]
Comparison then gives $0 \leq nu_t^{f/n}\le(n+1)u_t^{f/(n+1)}\le v_t^f$.
	Integrating over space and letting $t\to\infty$, we obtain
	\[
	0\le nU_\infty(1/n)
	\le(n+1)U_\infty(1/(n+1))
	\le V_\infty^f
	\qquad\text{a.s.}.
	\]
	Since $n$ ranges over a countable set, these inequalities hold
	simultaneously for all $n\ge1$ on a common event of probability one. Consequently, the limit
	\[
	H(f):=\lim_{n\to\infty}nU_\infty(1/n)
	\]
	exists and satisfies $0\le H(f)\le V_\infty^f$ almost surely.

For $\lambda>0$, define
$D_t^\lambda=V_t^f-U_t(\lambda)/\lambda\ge0$. This difference satisfies
\[
 D_t^\lambda=\frac1{2\lambda}\int_0^t\|u_s^{\lambda f}\|_2^2ds
 +\sqrt\kappa\int_0^t\!\int
             (v_s^f-u_s^{\lambda f}/\lambda)(x)\,W(ds,dx).
\]
For each integer $N\ge1$, let
\[
 \sigma_N=\inf\left\{t\ge0:\int_0^t\|v_s^f\|_2^2\,ds\ge N\right\},
 \qquad \inf\varnothing=\infty.
\]
Let $\mathcal N^\lambda$ denote the stochastic integral in the equation
for $D^\lambda$. Since $0\le v-u/\lambda\le v$ and $0\le q\le1$, we have
\begin{align*}
 \E\big[[\mathcal N^\lambda]_{t\wedge\sigma_N}\big]
 &=\kappa\E\int_0^{t\wedge\sigma_N}
                Q(v_s^f-u_s^{\lambda f}/\lambda)\,ds\\
 &\le\kappa\int_0^t\E(V_s^f)^2\,ds
 \le\kappa F^2\int_0^t e^{\kappa s}\,ds<\infty.
\end{align*}
 Consequently,
\begin{align*}
 \E D_{t\wedge\sigma_N}^\lambda
 &=\frac1{2\lambda}\E\int_0^{t\wedge\sigma_N}
                            \|u_s^{\lambda f}\|_2^2ds\\
 &\le\frac\lambda2\E\int_0^{t\wedge\sigma_N}\|v_s^f\|_2^2ds
 \le\frac{\lambda N}{2}.
\end{align*}
Set $D_\infty^\lambda=V_\infty^f-U_\infty(\lambda)/\lambda$.
Since $D^\lambda\ge0$, Fatou's lemma gives the explicit bound
\begin{align*}
 \E[D_\infty^\lambda\mathbf1_{\{\sigma_N=\infty\}}]
 &=\E\left[\liminf_{t\to\infty}
          D_{t\wedge\sigma_N}^\lambda
                      \mathbf1_{\{\sigma_N=\infty\}}\right]\\
 &\le\liminf_{t\to\infty}\E D_{t\wedge\sigma_N}^\lambda
 \le\frac{\lambda N}{2}.
\end{align*}
Taking $\lambda=1/n$ yields
\[
 0\le\E[(V_\infty^f-H(f))\mathbf1_{\{\sigma_N=\infty\}}]
 \le\liminf_{n\to\infty}\frac{N}{2n}=0.
\]
Thus $H(f)=V_\infty^f$ on each $\{\sigma_N=\infty\}$.
By Lemma~\ref{lem:linear-energy}, the union of these events has
probability one, which proves the asserted identity.

Finally, taking $\lambda=1/n$ in \eqref{eq:pi-laplace} gives
\[
 \int n(1-e^{-\mu(f)/n})\,\pi_\kappa(d\mu)
       =\E\big[n(1-e^{-U_\infty(1/n)})\big].
\]
The integrands on the left increase to $\mu(f)$. 
We have $nU_\infty(1/n) \uparrow V_\infty^f$ almost surely. Hence
\[
0\le U_\infty(1/n)
\le\frac{V_\infty^f}{n}\longrightarrow0
\qquad\text{a.s.}.
\]
It follows that
\[
n\bigl(1-e^{-U_\infty(1/n)}\bigr)
=nU_\infty(1/n) \int_0^1 e^{-\theta U_\infty(1/n)}\,d\theta
\longrightarrow V_\infty^f
\qquad\text{a.s.}.
\]
Moreover, since $0\le1-e^{-r}\le r$ for $r\ge0$, we have
\[
0\le n\bigl(1-e^{-U_\infty(1/n)}\bigr)
\le nU_\infty(1/n) \le V_\infty^f.
\]
Since $\E V_\infty^f\le F<\infty$, the dominated convergence theorem gives
\[
\lim_{n\to\infty}
\E\!\left[n\bigl(1-e^{-U_\infty(1/n)}\bigr)\right]
=\E V_\infty^f.
\]
On the other hand, $n(1-e^{-\mu(f)/n})\uparrow\mu(f)$, so
monotone convergence yields
\[
\lim_{n\to\infty}
\int n(1-e^{-\mu(f)/n})\,\pi_\kappa(d\mu)
=\int\mu(f)\,\pi_\kappa(d\mu).
\]
Combining these two limits 
proves \eqref{meanbridge}.
\end{proof}

\subsection{The terminal mean and the polymer disorder threshold}

This subsection proves that the terminal mean coefficient is either
zero or one and establishes a common disorder threshold for
integrable densities and point initial measures.
The terminal mean also determines the mean measure of the limiting law $\pi_\kappa$.

For a
signed $g\in L^1$, define $V_t^g:=V_t^{g^+}-V_t^{g^-}$ and its almost
sure terminal value by the corresponding difference. At every finite
time, \eqref{l2} gives
\[
 |V_t^g-V_t^h|
 =\left|\int(g-h)(x)Z_t^\kappa(x)\,dx\right|
 \le\int|g-h|(x)Z_t^\kappa(x)\,dx=V_t^{|g-h|}.
\]
Fatou's lemma therefore gives
\begin{equation}\label{eq:terminal-contraction}
 \E|V_\infty^g-V_\infty^h|
       \le\liminf_{t\to\infty}\E V_t^{|g-h|}
       =\|g-h\|_1.
\end{equation}
Thus $\Lambda_\kappa:g\mapsto V_\infty^g$ is a positive bounded linear
map from $L^1(dx)$ to $L^1(\Pp)$. Here linearity is understood in $L^1(\Pp)$.
 
 Let $s\ge0$ be deterministic, and let $G$ be an
$\mathcal F_s^W$-measurable $L^1_+$-valued random variable satisfying
$\E\|G\|_1<\infty$.  For $t\ge0$, let $V_t^{(s),G}$ denote the total mass at time $s+t$
of the linear equation started from the density $G$ at time $s$.

The following lemma gives the terminal mean formulas for deterministic
and random initial data.
\begin{lemma}\label{lem:terminal-mean}
	There exists a constant $c_\kappa\in[0,1]$ such that
	\begin{equation}\label{eq:terminal-mean-coefficient}
		\E V_\infty^g=c_\kappa\int g,
		\qquad g\in L^1(\mathbb R^d).
	\end{equation}
	Moreover,
 $(V_t^{(s),G})_{t\ge0}$ is a nonnegative martingale with respect
	to $(\mathcal F_{s+t}^W)_{t\ge0}$ and has a finite almost sure limit
	$V_\infty^{(s),G}$ satisfying
	\begin{equation}\label{eq:random-initial-mean}
		\E[V_\infty^{(s),G}\mid\mathcal F_s^W]
		=c_\kappa\|G\|_1.
	\end{equation}
\end{lemma}

\begin{proof}
	Define 
	\[
	\ell_\kappa(g):=\E\Lambda_\kappa(g)=\E V_\infty^g,
	\qquad g\in L^1(\mathbb R^d).
	\]
	The estimate \eqref{eq:terminal-contraction} gives
	$|\ell_\kappa(g)|\le\|g\|_1$.
	By stationarity of the environment, $\ell_\kappa$ is invariant
	under spatial translations, so
	\[
	\ell_\kappa(g(\cdot-y))=\ell_\kappa(g),
	\qquad y\in\mathbb R^d.
	\]
	Choose $\phi\in C_c^{\infty,+}$ with $\int\phi=1$.
Linearity and translation invariance give
	\begin{align*}
		\ell_\kappa(g*\phi)
		&=\int\phi(y)\ell_\kappa(g(\cdot-y))\,dy
		=\ell_\kappa(g),\\
		\ell_\kappa(g*\phi)
		&=\int g(y)\ell_\kappa(\phi(\cdot-y))\,dy
		=\ell_\kappa(\phi)\int g.
	\end{align*}
	Hence \eqref{eq:terminal-mean-coefficient} holds with
	$c_\kappa:=\ell_\kappa(\phi)$. Positivity and the contraction bound
	imply $0\le c_\kappa\le\|\phi\|_1=1$.

	Let $w_t$ denote the density at time $s+t$ of the linear evolution started from $G$ at time $s$. This density is
	$\mathcal F_{s+t}^W$-measurable. The kernel
	$\mathcal K_{s+t,s+t+r}^\kappa$ depends only on future noise,
	is independent of $\mathcal F_{s+t}^W$, and satisfies
	\[
	\E\left[\int\mathcal K_{s+t,s+t+r}^\kappa(x,y)\,dy\right]=1.
	\]
	Thus, by the flow property and conditional Tonelli's theorem,
	\begin{align*}
		\E[V_{t+r}^{(s),G}\mid\mathcal F_{s+t}^W]
		&=\int w_t(x)\,
		\E\left[\int\mathcal K_{s+t,s+t+r}^\kappa(x,y)\,dy\right]dx\\
		&=\int w_t(x)\,dx
		=V_t^{(s),G},
		\qquad t\ge0,\ r>0.
	\end{align*}
	
Moreover,
	\[
	\E[V_t^{(s),G}\mid\mathcal F_s^W]=\|G\|_1,
	\qquad
	\E V_t^{(s),G}=\E\|G\|_1<\infty.
	\]
	The nonnegative martingale $V^{(s),G}$ therefore has a finite almost sure limit $V_\infty^{(s),G}$, and Fatou's lemma gives
	\[
	\E V_\infty^{(s),G}\le\E\|G\|_1<\infty.
	\]

	Since $L^1(\mathbb R^d)$ is separable, we may choose
	$\mathcal F_s^W$-measurable simple random variables
	\[
	G_j=\sum_{i=1}^{N_j}\mathbf1_{A_{j,i}}g_{j,i},
	\qquad g_{j,i}\in L^1_+,
	\]
	where $(A_{j,i})_{i=1}^{N_j}$ is a measurable partition. We choose the approximations so that
	$
	\varepsilon_j:=\E\|G_j-G\|_1\longrightarrow0.
	$
	
	Independence and stationarity of the future noise, together with
	\eqref{eq:terminal-mean-coefficient}, yield
	\begin{align*}
		\E[V_\infty^{(s),G_j}\mid\mathcal F_s^W]
		&=\sum_{i=1}^{N_j}\mathbf1_{A_{j,i}}
		\E V_\infty^{g_{j,i}}
		=c_\kappa\sum_{i=1}^{N_j}\mathbf1_{A_{j,i}}\|g_{j,i}\|_1
		=c_\kappa\|G_j\|_1.
	\end{align*}
	
	Note that, for every finite $t$,
	\[
	|V_t^{(s),G_j}-V_t^{(s),G}|
	\le V_t^{(s),|G_j-G|}.
	\]
	We obtain
	\[
	\E|V_t^{(s),G_j}-V_t^{(s),G}|
	\le\E\|G_j-G\|_1=\varepsilon_j.
	\]
Fatou's
	lemma then gives
	\[
	\E|V_\infty^{(s),G_j}-V_\infty^{(s),G}|
	\le\liminf_{t\to\infty}
	\E|V_t^{(s),G_j}-V_t^{(s),G}|
	\le\varepsilon_j.
	\]
	Consequently, the $L^1$ contraction property of conditional
	expectation implies
	\begin{align*}
		&\left\|
		\E[V_\infty^{(s),G}\mid\mathcal F_s^W]
		-c_\kappa\|G\|_1
		\right\|_{L^1(\Pp)}\\
		&\quad\le
		\E|V_\infty^{(s),G}-V_\infty^{(s),G_j}|
		+c_\kappa\E\big|\|G_j\|_1-\|G\|_1\big|\\
		&\quad\le(1+c_\kappa)\varepsilon_j
		\longrightarrow0.
	\end{align*}
	This proves \eqref{eq:random-initial-mean}.
\end{proof}

Fix $\kappa\ge0$. For a finite deterministic nonnegative measure
$\nu$, define
\[
Z_t^\nu:=\int_{\mathbb R^d}Z_t^\kappa(x)\,\nu(dx),
\qquad t\ge0.
\]
This is the total mass of the linear evolution in \eqref{eq:dual-equations} started from $\nu$.

We next show that $c_\kappa$ is either $0$ or $1$ and derive the
corresponding conclusions for the limiting law and the total mass
of the linear evolution from finite initial measures.

\begin{proposition}
\label{prop:mean-dichotomy}
For every $\kappa\ge0$,  
 $c_\kappa\in\{0,1\}$. The limiting law
$\pi_\kappa$ equals $\delta_0$ when $c_\kappa=0$ and has mean measure $m$ when
$c_\kappa=1$. For every finite deterministic nonnegative measure $\nu$,
the total mass $Z_t^\nu$ satisfies
\[
 \E Z_\infty^\nu=c_\kappa\nu(1).
\]
The martingale is uniformly integrable when $c_\kappa=1$, and its terminal value is
zero almost surely when $c_\kappa=0$.
\end{proposition}
\begin{proof}
Take $f\in C_c^{\infty,+}$ with $m(f)>0$. By \eqref{eq:linear-mass-moments}, the $L^1$ norm of $v_s^f$ is integrable. The flow property and
\eqref{eq:random-initial-mean} give
\[
 \E[V_\infty^f\mid\mathcal F_s^W]=c_\kappa V_s^f.
\]
 Letting $s=n\to\infty$ and applying L\'evy's upward theorem, we obtain
 $(1-c_\kappa)V_\infty^f=0$ almost surely. Taking expectations
and using \eqref{eq:terminal-mean-coefficient} gives
\[
 0=(1-c_\kappa)\E V_\infty^f
   =c_\kappa(1-c_\kappa)m(f).
\]
Since $m(f)>0$ and $c_\kappa\in[0,1]$, this proves
$c_\kappa\in\{0,1\}$.

By \eqref{meanbridge} and \eqref{eq:terminal-mean-coefficient}, the mean
measure of $\pi_\kappa$ agrees with $c_\kappa m$ on $C_c^{\infty,+}$.
A smooth majorant of the indicator of each compact set gives local
finiteness, so uniqueness of Radon measures gives equality on all Borel sets. When $c_\kappa=1$, the limiting law has mean measure $m$ and the population is persistent.
When $c_\kappa=0$, we have $\pi_\kappa=\delta_0$, so $X_t(f)\to0$ in probability for every $f\in C_c^+(\R^d)$.

At time $1$, the linear evolution in \eqref{eq:dual-equations} started from $\nu$
has the density
\[
G(y):
=\int \mathcal K_{0,1}^\kappa(x,y)\,\nu(dx).
\]
This density is nonnegative and $\mathcal F_1^W$-measurable, with
\[
\|G\|_1=Z_1^\nu,
\qquad
\E\|G\|_1=\nu(\mathbb R^d)<\infty.
\]
Thus $G\in L^1_+$ almost surely. By the flow property and
\eqref{eq:random-initial-mean}, we have
\[
\E[Z_\infty^\nu\mid\mathcal F_1^W]
=c_\kappa\|G\|_1
=c_\kappa Z_1^\nu.
\]

Taking expectations gives
\[
 \E Z_\infty^\nu=c_\kappa\E Z_1^\nu=c_\kappa\nu(1).
\]
When $c_\kappa=1$, we have $\E Z_t^\nu=\E Z_\infty^\nu=\nu(1)$, so
\[
\E|Z_t^\nu-Z_\infty^\nu|
=\E Z_t^\nu+\E Z_\infty^\nu
-2\E\min(Z_t^\nu,Z_\infty^\nu)
\longrightarrow0.
\]
Thus $Z_t^\nu\to Z_\infty^\nu$ in $L^1(\Pp)$, and the martingale
$(Z_t^\nu)_{t\ge0}$ is uniformly integrable. When $c_\kappa=0$,
$\E Z_\infty^\nu=0$ and $Z_\infty^\nu\ge0$ force
$Z_\infty^\nu=0$ almost surely. In particular, all statements apply to
$\nu=\delta_x$.
\end{proof}
Define
\begin{equation}\label{eq:critical-parameter}
	\kappa_c:=\sup\{\kappa\ge0:c_\kappa=1\}.
		\end{equation}
The following proposition establishes a finite threshold between the two regimes.
\begin{proposition}\label{prop:critical-parameter}
	The map $\kappa\mapsto c_\kappa$ is nonincreasing. Moreover, $\kappa_c<\infty$, and
	\[
	c_\kappa=1\quad\text{for }0\le\kappa<\kappa_c,
	\qquad
	c_\kappa=0\quad\text{for }\kappa>\kappa_c.
	\]
\end{proposition}

\begin{proof}
	Let $0\le\kappa<\kappa'$, and let $W_1,W_2$ be independent
	copies of the environment. Construct $Z_t^\kappa$ with driving noise $\sqrt\kappa W_1$ and $Z_t^{\kappa'}$ with driving noise
	$\sqrt\kappa W_1+\sqrt{\kappa'-\kappa}W_2$. 
	Taking expectations over $W_2$ in the Feynman--Kac representation
	of $Z_t^{\kappa'}$ gives
	\[
	\E[Z_t^{\kappa'}\mid\sigma(W_1)]=Z_t^\kappa.
	\]
	Therefore, we have
	\[
	\E[Z_\infty^{\kappa'}\mid\sigma(W_1)]
	\le\liminf_{n\to\infty}
	\E[Z_n^{\kappa'}\mid\sigma(W_1)]
	=Z_\infty^\kappa.
	\]
	Taking expectations and using Proposition~\ref{prop:mean-dichotomy},
	we obtain $c_{\kappa'}\le c_\kappa$.
	
Since $q\in C_b^\infty$ and $q(0)=1$,
we may apply \cite[Theorem~1.5]{CRZ} to obtain a finite constant $N_0$ such that local extinction
holds whenever $\kappa\ge N_0$.
For such $\kappa$, we have $\pi_\kappa=\delta_0$, and
therefore $c_\kappa=0$. Since $\{\kappa\ge0:c_\kappa=1\}\subset[0,N_0)$, we obtain
\[ 
\kappa_c\le N_0<\infty.
\]

If $0\le\kappa<\kappa_c$, the definition of the supremum ensures
that there exists $\kappa'>\kappa$ with $c_{\kappa'}=1$.
Monotonicity then gives $c_\kappa\ge c_{\kappa'}=1$,
so $c_\kappa=1$.
If $\kappa>\kappa_c$, then $c_\kappa\ne1$ by the definition of
$\kappa_c$. Hence $c_\kappa=0$ by Proposition~\ref{prop:mean-dichotomy}.
\end{proof}

\section{Exponential decay of fractional moments above the critical value}\label{sec:sharpness}
We prove that if $\kappa>b_*$ and $Z_t\to0$ almost surely, then
$\E[Z_t^\theta]$ decays exponentially in $t$ for every $\theta\in(0,1)$.
These hypotheses hold for every $\kappa>\kappa_c$.
Throughout this section, fix $\kappa>0$, set $\beta=\sqrt\kappa$,
and write $\mathcal F_t=\mathcal F_t^W$.
Denote the two heat semigroups by
\[
(P_t^Bf)(x)=E_x^B[f(B_t)],
\qquad
(P_t^Df)(x)=E_x^D[f(D_t)].
\]

For $\nu\in\Mp$,
let $M_t^\nu(A)$ denote the mass in a Borel set $A$ at time $t$
for the linear equation in \eqref{eq:dual-equations}
 started from $\nu$. Its Feynman--Kac representation is
\begin{equation}\label{eq:sharpness-forward-measure}
	M_t^\nu(A)
	=\int\nu(dx)E_x^B\left[
	e^{\beta H_t(B)-\kappa t/2}\mathbf1_{\{B_t\in A\}}\right],
	\qquad A\in\mathcal B(\R^d).
\end{equation}
Its total mass is $M_t^\nu(\R^d)=Z_t^\nu$. Moreover, $Z^\nu$ is a continuous positive martingale with mean one.
For every $T<\infty$,
\begin{equation}\label{eq:sharpness-replica-finite}
	\sup_{t\le T}\E[(Z_t^\nu)^2]
	=\sup_{t\le T}\iint\nu(dx)\nu(dy)E_{x,y}^{B,B'}
	e^{\kappa\int_0^tq(B_s-B'_s)ds}\le e^{\kappa T}.
\end{equation}
With $\kappa$ fixed, we recall the notation
\[
Z_t(x):=Z_t^\kappa(x)=Z_t^{\delta_x},
\qquad
Z_t:=Z_t(0)=Z_t^\kappa.
\]

We record the elementary evolution properties needed below. For bounded
$C^2$ test functions with bounded derivatives, the forward measure satisfies
\begin{equation}\label{eq:sharpness-weak-flow}
 d\langle M_t^\nu,\varphi\rangle
 =\tfrac12\langle M_t^\nu,\Delta\varphi\rangle\,dt
 +\beta\int\varphi(x)M_t^\nu(dx)W(dt,x),
\end{equation}
and the corresponding mild equation is
\[
\langle M_t^\nu,\varphi\rangle
=\langle\nu,P_t^B\varphi\rangle
+\beta\int_0^t\int(P_{t-s}^B\varphi)(x)M_s^\nu(dx)W(ds,x),
\]
which has a unique nonnegative solution.
For any almost surely finite stopping time $\tau$,
the shifted noise
\[
\widetilde W_t(x):=W_{\tau+t}(x)-W_\tau(x),\qquad t\ge0,
\]
is independent of $\mathcal F_\tau$ and has the same law as $W$.
The continuation $(M_{\tau+t}^\nu)_{t\ge0}$ coincides with the
linear evolution started from $M_\tau^\nu$ and driven by
$\widetilde W$.

Recall the definitions of $b_*$ and the operator
$T$ in \eqref{eq:second-moment-threshold-intro}. We have the following result. 

\begin{theorem}\label{thm:sharp}
Suppose that $\kappa>b_*$ and $Z_t\to0$ almost surely. Then, for every $\theta\in(0,1)$,
there are constants $C,c\in(0,\infty)$ such that
\[
 \E[Z_t^\theta]\le Ce^{-ct},\qquad t\ge0.
\]
Moreover,
\[
 \limsup_{t\to\infty}\frac1t\E\log Z_t<0.
\]
In particular, both conclusions hold for every $\kappa>\kappa_c$.
\end{theorem}

The proof follows the excursion and size-bias strategy of
Junk and Lacoin~\cite{JL}. We construct an auxiliary kernel and show how polynomial
fractional-moment decay implies exponential decay in the present
spatially correlated Gaussian environment.

\begin{lemma}\label{lem:spectral}
The operator $T$ is positive and compact on $C_0(\R^d)$, with
$0<r(T)<\infty$. Moreover, $b_*\le\kappa_c$.
For every $\kappa>b_*$, there exist a real-valued even function
$h\in C_0(\R^d)\cap C_b^2(\R^d)$ and a constant $a>0$ such that
\begin{equation}\label{eq:sharpness-auxiliary-inequality}
	(\Delta+\kappa q)h\ge aq\quad\text{on }\R^d.
\end{equation}
\end{lemma}
\begin{proof}
Choose smooth radial cutoffs $0\le\chi_R\le1$, equal to one on $B_R$
and zero on $B_{2R}^c$, increasing to one as $R\to\infty$. Define
$q_R=q\chi_R$ and $q_R^{\rm tail}=q(1-\chi_R)$. The Newton kernel of
$-\Delta$ is a positive constant times $|x|^{2-d}$. By radial monotonicity
and scaling, $\sup_x\int_{B_r}|x-y|^{2-d}\,dy\le Cr^2$. The annulus
$2^jR\le |y|<2^{j+1}R$ therefore contributes at most
$C(2^jR)^{2-\alpha}$ for $R\ge1$. Summing over $j\ge0$ gives
\begin{equation}\label{eq:sharpness-green-tail}
 \|Gq_R^{\rm tail}\|_\infty\le CR^{2-\alpha}\longrightarrow0.
\end{equation}
Note that
$Gq_R(x)=O_R(|x|^{2-d})$ as $|x|\to\infty$.
Thus $Gq_R\in C_0(\R^d)$. By \eqref{eq:sharpness-green-tail}, we have
\[
\|Gq-Gq_R\|_\infty
=\|Gq_R^{\rm tail}\|_\infty\longrightarrow0.
\]
Hence $Gq$, being a uniform limit of functions in $C_0(\R^d)$,
also belongs to $C_0(\R^d)$.

Define $T_Rf:=G(q_Rf)$. For every
$f\in C_0(\R^d)$ with $\|f\|_\infty\le1$, we have
\begin{align*}
 |T_Rf(x)|&\le C_R(|x|-2R)^{2-d},&& |x|>2R+1,\\
 |T_Rf(x+z)-T_Rf(x)|
 &\le C\eta^2+C_R\eta^{1-d}|z|,&& |z|<\eta/2,
\end{align*}
where $0<\eta\le1$ and the constants are independent of
$f,x,z$ and $\eta$. The first estimate shows that $T_Rf$ vanishes at infinity uniformly
over the unit ball of $C_0(\R^d)$.
In the second estimate, choosing $\eta$ small and then $|z|$ small
proves equicontinuity of the same family. The Arzel\`a--Ascoli theorem now gives compactness of $T_R$, and
\[
 \|T-T_R\|\le\|Gq_R^{\rm tail}\|_\infty\longrightarrow0
\]
proves compactness of $T$.

Choose $0\le g\in C_c(\R^d)$, $g\not\equiv0$,
and let $c=\min_{x\in \operatorname{supp}g}Tg(x)/\|g\|_\infty>0$. Positivity then yields
\[
 Tg\ge cg,\qquad T^ng\ge c^ng,
 \qquad 0 <c\le\lim_{n\to\infty}\|T^n\|^{1/n}=r(T)\le\|Gq\|_\infty.
\]
This proves the first assertion of the lemma.

For $j\ge1$, set $\Gamma_j:=T^{j-1}(Gq)$. Tonelli's theorem and the Markov property give
\[
 \Gamma_j(x)=\int_{0<s_1<\cdots<s_j<\infty}
 E_x^D\!\left[\prod_{i=1}^j q(D_{s_i})\right]ds_1\cdots ds_j.
\]
Therefore,
\begin{align*}
 E_x^D e^{b\int_0^\infty q(D_s)ds}
 &=1+\sum_{j\ge1}b^j\Gamma_j(x)
  =1+\sum_{j\ge1}b^jT^{j-1}Gq(x).
\end{align*}
For $b \in (0, b_*)$, choose $\varrho\in(br(T),1)$. The spectral radius formula gives
$\|(bT)^j\|\le C\varrho^j$, and hence
\[
 \sup_x E_x^D e^{b\int_0^\infty q(D_s)ds}
 \le1+\frac{Cb\|Gq\|_\infty}{1-\varrho}<\infty.
\]
The second-moment identity and the preceding estimate give
\[
\sup_{t\ge0}\E[(Z_t^b)^2]
\le E_0^D\exp\left\{b\int_0^\infty q(D_s)\,ds\right\}<\infty.
\]
Hence $(Z_t^b)_{t\ge0}$ is uniformly integrable.
Applying Proposition~\ref{prop:mean-dichotomy} with parameter
$\kappa=b$ and initial measure $\nu=\delta_0$, we obtain
\[
c_b=\E Z_\infty^b=\lim_{t\to\infty}\E Z_t^b=1.
\]
Since $b$ was arbitrary in $[0,b_*)$, the definition
\eqref{eq:critical-parameter} yields $\kappa_c\ge b_*$.

Fix $\kappa>b_*$. By the compactness and positivity of $T$, the
Krein--Rutman theorem gives a nonzero $\psi\in C_0(\R^d)$ with
$\psi\ge0$ and $\psi=b_*G(q\psi)$. The Green kernel and $q$ are
strictly positive, so $\psi(x)>0$ for every $x$. Since $q$ is even,
we may replace $\psi(x)$ by $(\psi(x)+\psi(-x))/2$.
Put $U=Gq$. Then $U\in C_0(\R^d)$ is even and
\[
 -\Delta U=q,\qquad -\Delta\psi=b_*q\psi
\]
in the sense of distributions. Both $U$ and $\psi$ belong to
$C_b^2(\R^d)$. Indeed, the functions and their right-hand sides
are bounded, so interior $W^{2,p}$ estimates on balls of fixed radii,
with $p>d$, give uniform local $C^{1,\gamma}$ bounds for some
$\gamma>0$. Since $q\in C_b^1$, the right-hand sides now have
uniform local $C^\gamma$ bounds, and interior Schauder estimates
give uniform bounds on the second derivatives.

Choose $R$ such that $\kappa U\le1/2$ on $B_R^c$ and put
$m_R=\min_{\overline B_R}\psi>0$. Take $A>0$ such that
$A(\kappa-b_*)m_R\ge\kappa\|U\|_\infty$, and set $h=A\psi-U$.
Then
\[
 (\Delta+\kappa q)h
 =q\bigl[A(\kappa-b_*)\psi+1-\kappa U\bigr].
\]
The bracket is at least one on $\overline B_R$ and at least
$1/2$ on $B_R^c$. Thus $h\in C_0(\R^d)\cap C_b^2(\R^d)$ is even
and satisfies \eqref{eq:sharpness-auxiliary-inequality} with $a=1/2$.
\end{proof}

\begin{lemma}[Comparison on a unit cube]\label{lem:derivatives}
Let $Q=[0,1]^d$. For every $0<\theta\le1$, there is a constant
$C_\theta<\infty$ such that
\begin{equation}\label{eq:sharpness-derivative-sup}
 \E\!\left[\left(\sup_{x\in Q} Z_t(x)\right)^\theta\right]
 \le C_\theta\E[Z_t(0)^\theta],\qquad t\ge1.
\end{equation}
There is also a constant $C<\infty$ such that
\begin{equation}\label{eq:sharpness-kernel-tail}
 \E\sup_{x\in Q}\mathcal K_{0,t}^\kappa(x,y)
 \le C p_{t+1}(y),\qquad t\ge1,\quad y\in\R^d .
\end{equation}
\end{lemma}

\begin{proof}
Write $K_{s,t}=\mathcal K_{s,t}^\kappa$. For $T>0$, let $b$ be a
Brownian bridge from $0$ to $0$ of duration $T$, independent of $W$.
Conditioning the reference Brownian path on its endpoint gives
\[
 \begin{aligned}
 K_{0,T}(x,y)&=p_T(y-x)E^{\rm br}
       \exp\{\beta H_{x,y}^b-\kappa T/2\},\\
 H_{x,y}^b&=\int_0^T W\bigl(ds,(1-s/T)x+(s/T)y+b_s\bigr).
 \end{aligned}
\]
Conditionally on $b$, the variables $H_{x,y}^b$ form a centered Gaussian field
whose covariance is independent of $b$. Since $1-q(z)\le C|z|^2$,
the standard deviation of each increment is at most
$C_T(|x-x'|+|y-y'|)$. Kolmogorov's continuity theorem gives a
continuous version on compact sets, obtained measurably by extension
from rational endpoints. Fernique's theorem gives finite exponential
moments of its supremum on each compact set, uniformly in $b$.
Fubini's theorem and dominated convergence therefore give a jointly measurable
version of the bridge average that is continuous in the endpoints
almost surely for each fixed $T$. For $Z_T(x)$, the same construction applies to the field
$x\mapsto H_T(x+B)$, whose conditional covariance given $B$ is $Tq(x-x')$.
Thus spatial suprema over $Q$ may be taken over a countable dense subset.

Set $\Gamma(y)=\sup_{x\in Q}K_{0,1}(x,y)$. For $T=1$ and fixed $y$,
the conditional covariance of $H_{x,y}^b$ and $H_{x',y}^b$ equals
$\int_0^1q((1-s)(x-x'))\,ds$, independently of both $b$ and $y$.
Taking the supremum inside the bridge expectation and using the
preceding exponential moment bound, we obtain
\begin{equation}\label{eq:unit-bridge-envelope}
 \E\Gamma(y)\le C\sup_{x\in Q}p_1(y-x)\le C p_2(y).
\end{equation}
For $t>1$, kernel composition, independence of the noise on the two time intervals,
and $\E K_{1,t}(z,y)=p_{t-1}(y-z)$ now give
\[
 \E\sup_{x\in Q}K_{0,t}(x,y)
 \le\int\E\Gamma(z)\,p_{t-1}(y-z)\,dz\le C p_{t+1}(y).
\]
The case $t=1$ follows from \eqref{eq:unit-bridge-envelope}.

For $t>1$, put $V(y)=\int K_{1,t}(y,z)\,dz$, and for $t=1$, put
$V\equiv1$. Conditional on $V$, concavity and
\eqref{eq:unit-bridge-envelope} give
\[
 \E\!\left[\left(\sup_{x\in Q}Z_t(x)\right)^\theta\mid V\right]
 \le\left(\int\E\Gamma(y)V(y)\,dy\right)^\theta
 \le C^\theta m^\theta,\qquad m=(P_2^B V)(0).
\]
Here $m<\infty$ almost surely since $\E m=1$. Let $K'$ be the kernel
of an independent environment on a time interval of length two, and
set $Y=\int K'_{0,2}(0,y)V(y)\,dy$. Then $\E[Y\mid V]=m$.
For two independent reference paths $B,B'$ started at $0$, integration
over the independent environment gives, conditionally on $V$,
\[
 \E[Y^2\mid V]
 = E_{0,0}^{B,B'}\!\left[V(B_2)V(B'_2)
       e^{\kappa\int_0^2q(B_s-B'_s)\,ds}\right]\le e^{2\kappa}m^2.
\]
The conditional Paley--Zygmund inequality thus gives
$\E[Y^\theta\mid V]\ge 2^{-\theta-2}e^{-2\kappa}m^\theta$
on $\{m>0\}$, and the same bound is immediate on $\{m=0\}$.
Concatenating the independent time intervals gives
$Y\stackrel{\rm law}{=}Z_{t+1}(0)$. Since $Z^\theta$ is a
supermartingale, $\E Y^\theta\le\E Z_t(0)^\theta$.
Combining the last two bounds proves \eqref{eq:sharpness-derivative-sup}.
\end{proof}

We next use the auxiliary kernel from Lemma~\ref{lem:spectral}
to define a bounded functional of the normalized mass distribution. Recall the definition of $M_t^\nu$ in \eqref{eq:sharpness-forward-measure}.

\begin{lemma}
	\label{lem:drift}
	Fix $\kappa>b_*$, and let $h,a$ be given by
	Lemma~\ref{lem:spectral}. For $\nu\in\Mp$, set
	\[
	\mu_t=M_t^\nu/Z_t^\nu,\qquad
		I^g_t=\iint g(x-y)\mu_t(dx)\mu_t(dy),\quad g\in\{h,q\}.
	\]
	Then $dI^h_t=b_t\,dt+dN_t$, where $N$ is a continuous local
	martingale. Moreover, there exist a constant $C>0$ and a deterministic
	function $e(\delta)\to0$ as $\delta\downarrow0$, both independent
	of $\nu$, such that whenever $I^q_t\le\delta$,
	\[
	b_t\ge(a-Ce(\delta))I^q_t,\qquad
	d\langle N\rangle_t\le Ce(\delta)^2I^q_t\,dt.
	\]
\end{lemma}

\begin{proof}
For a probability measure $\mu$, write
$I^q(\mu)=\langle\mu,q*\mu\rangle$ and set
$c_q=\inf_{|z|\le2}q(z)>0$. Since $q\ge0$, we have
\begin{equation}\label{md}
 I^q(\mu)\ge\iint_{B(x,1)^2}q(y-z)\mu(dy)\mu(dz)
       \ge c_q\mu(B(x,1))^2.
\end{equation}
By \eqref{md}, we have
\begin{equation}\label{uy}
 \|h*\mu\|_\infty
 \le C_L\|h\|_\infty I^q(\mu)^{1/2}
             +\sup_{|y|>L}|h(y)|.
\end{equation}
Define
\[
 e(\delta)=\sup\{\|h*\mu\|_\infty:
                  \mu\in\Mp,\ I^q(\mu)\le\delta\}.
\]
For fixed $L$, first let $\delta\downarrow0$ in \eqref{uy} and then let
$L\to\infty$. Since $h\in C_0(\R^d)$, this proves
$e(\delta)\downarrow0$. Also, $|\langle\mu,h*\mu\rangle|\le e(\delta)$
whenever $I^q(\mu)\le\delta$.

Define
\[
V_t:=\iint h(x-y)\,M^\nu_t(dx)M^\nu_t(dy).
\]
Since $\mu_t=M^\nu_t/Z^\nu_t$, we have $I_t^h=V_t/(Z^\nu_t)^2$.

Taking the constant test function $1$ in
\eqref{eq:sharpness-weak-flow} gives
\[
dZ_t^\nu=\beta\int M_t^\nu(dx)\,W(dt,x).
\]
To compute $dV_t$, first apply the product rule in
\eqref{eq:sharpness-weak-flow} to test functions $\phi(x)\xi(y)$.
Stop when $Z_t^\nu$ first leaves $(1/m,m)$, $m>1$. For fixed $R$,
the kernel $\chi_R(x)\chi_R(y)h(x-y)$ is a $C^2$ limit of finite sums
of smooth compactly supported tensors, by convolution with a product
mollifier followed by Riemann sums. As $R\to\infty$, the kernels and
their derivatives through order two are uniformly bounded and converge
pointwise to those of $h(x-y)$. On each stopped finite time interval, the
absolute drift density is bounded by $Cm^2$ and the bracket density
by $C\kappa m^4$. Dominated convergence and the It\^o isometry justify
both limits. Since $Z^\nu$ is continuous and strictly positive at all
finite times, these stopping times increase to infinity. Using the
evenness of $h$ to combine the two noise terms gives
\[
\begin{aligned}
	dV_t
	={}&\iint (\Delta h+\kappa qh)(x-y)\,
	M_t^\nu(dx)M_t^\nu(dy)\,dt
	+2\beta\int (h*M_t^\nu)(x)\,M_t^\nu(dx)\,W(dt,x).
\end{aligned}
\]
Thus
the quadratic variation and covariation satisfy
\begin{align*}
	d\langle Z^\nu\rangle_t
	&=\kappa\iint q(x-y)\,M_t^\nu(dx)M_t^\nu(dy)\,dt
	=\kappa (Z_t^\nu)^2 I_t^q\,dt,\\
	d\langle V,Z^\nu\rangle_t
	&=2\kappa\int (h*M_t^\nu)(x)(q*M_t^\nu)(x)\,M_t^\nu(dx)\,dt.
\end{align*}

Applying It\^o's formula
to $I_t^h=V_t/(Z_t^\nu)^2$ now gives
\[
dI_t^h=b_t\,dt+dN_t,
\]
where
\begin{align*}
	b_t
	={}&\iint (\Delta h+\kappa qh)(x-y)\,\mu_t(dx)\mu_t(dy)
	+3\kappa I_t^q I_t^h\\
	&-4\kappa\int (h*\mu_t)(x)(q*\mu_t)(x)\,\mu_t(dx),\\
	dN_t
	={}&2\beta\int \bigl((h*\mu_t)(x)-I_t^h\bigr)\,
	\mu_t(dx)\,W(dt,x).
\end{align*}

If $I^q_t\le\delta$, then $|I^h_t|\le e(\delta)$ and
$\|h*\mu_t\|_\infty\le e(\delta)$. Using $q*\mu_t\ge0$ and
\eqref{eq:sharpness-auxiliary-inequality}, we therefore have
\[
 b_t\ge aI^q_t-3\kappa e(\delta)I^q_t
                   -4\kappa e(\delta)I^q_t
       =(a-7\kappa e(\delta))I^q_t.
\]
For $g_t(x)=(h*\mu_t)(x)-I^h_t$, the bound $|g_t(x)|\le2e(\delta)$ gives
\begin{align*}
 d\langle N\rangle_t
 &=4\kappa\iint g_t(x)g_t(y)q(x-y)\mu_t(dx)\mu_t(dy)\,dt\\
 &\le16\kappa e(\delta)^2 I^q_t\,dt.
\end{align*}
All estimates are uniform in $\nu$.
\end{proof}

The following lemma estimates the probability of a large accumulation of mass in a ball of fixed radius.

\begin{lemma}\label{lem:excursion}
Assume that $\kappa>b_*$ and $Z_t\to0$ almost surely. There exist
$r,c,v_0>0$ such that, for every $\nu\in\Mp$ and $v\ge v_0$,
\begin{equation}\label{eq:1110a}
 \Pp\left(\sup_{t\ge0,\ z\in\Z^d}M_t^\nu(B(z,r))\ge v\right)\ge c/v.
\end{equation}
\end{lemma}

\begin{proof}
	Fix $\nu\in\Mp$. Take $h$ and $a$ from Lemma~\ref{lem:spectral}
	at this value of $\kappa$. By
	Proposition~\ref{prop:mean-dichotomy}, $Z_t^\nu\to0$ almost surely.
	
	For $u>1$, let $\tau_u:=\inf\{t\ge0:Z_t^\nu\ge u\}$.
	The stopped martingale $Z_{t\wedge\tau_u}^\nu$ is bounded by $u$
	and converges to $u\mathbf1_{\{\tau_u<\infty\}}$. Since
	$Z_0^\nu=1$, bounded convergence gives
	\begin{equation}\label{eq:sharpness-exact-hit}
		\Pp(\tau_u<\infty)=\frac1u.
	\end{equation}
	
	Write $\mu_t=M_t^\nu/Z_t^\nu$ and $H=\|h\|_\infty$, so that
	$|I_t^h|\le H$. By Lemma~\ref{lem:drift}, we may choose
	$\delta\in(0,1)$ such that
	\[
	b_t\ge \frac a2 I_t^q
	\qquad\text{whenever } I_t^q\le\delta.
	\]
	Let $\Phi(z)=z-1-\log z$ for $z>0$. Since $\Phi(z)\to\infty$
	as $z\downarrow0$, choose $\eta\in(0,1)$ such that
	$(a/\kappa)\Phi(\eta)>4H$.
	
	By \eqref{eq:sharpness-weak-flow}, $M_t^\nu$ is vaguely continuous. Its total mass $Z_t^\nu$ is continuous and strictly positive, so $\mu_t=M_t^\nu/Z_t^\nu$ is weakly continuous. Since $(x,y)\mapsto q(x-y)$ is bounded and continuous, $I_t^q=\iint q(x-y)\,\mu_t(dx)\mu_t(dy)$ has continuous paths. On $\{\tau_u<\infty\}$, let
	\[
	\sigma:=\inf\{t\ge\tau_u:
	I_t^q\ge\delta \ \text{or}\ Z_t^\nu\le\eta u\}.
	\]
	This time is finite because $Z_t^\nu\to0$. 
	
	Define
	\[
	Y_t:=H-I_t^h+\frac a\kappa\Phi(Z_t^\nu/u).
	\]
	Since $\Phi''(z)=z^{-2}$ and
	$d\langle Z^\nu\rangle_t=\kappa(Z_t^\nu)^2I_t^q\,dt$,
	It\^o's formula gives the drift of $Y_t$ on $[\tau_u,\sigma)$ as
	\[
	-b_t+\frac a\kappa\,
	\frac{\Phi''(Z_t^\nu/u)}{2u^2}\,
	\kappa(Z_t^\nu)^2I_t^q
	=-b_t+\frac a2 I_t^q\le0.
	\]
	Set $\sigma=\infty$ on $\{\tau_u=\infty\}$. For $m\ge1$, let $\tau_m=\tau_u\wedge m$ and $A_m=\{\tau_u\le m\}\in\mathcal F_{\tau_m}$. In the filtration $(\mathcal F_{\tau_m+s})_{s\ge0}$, the process $\mathbf1_{A_m}Y_{(\tau_m+s)\wedge\sigma}$ is a nonnegative local supermartingale with initial value at most $2H$, hence a supermartingale. Since $\sigma<\infty$ on $A_m$, conditional Fatou's lemma gives the following inequality on $A_m$, and then on $\{\tau_u<\infty\}$ as $m\to\infty$.
	\[
	\E[Y_\sigma\mid\mathcal F_{\tau_u}]
	\le Y_{\tau_u}=H-I_{\tau_u}^h\le2H
	\qquad\text{on }\{\tau_u<\infty\}.
	\]
	On $\{Z_\sigma^\nu=\eta u\}$, we have
	$Y_\sigma\ge(a/\kappa)\Phi(\eta)>4H$. Consequently,
	\[
	\Pp\bigl(I_\sigma^q\ge\delta,\ Z_\sigma^\nu>\eta u
	\mid\mathcal F_{\tau_u}\bigr)\ge\tfrac12
	\qquad\text{on }\{\tau_u<\infty\}.
	\]
	
	Finally, choose $R>0$ such that $\sup_{|x|\ge R}q(x)\le\delta/2$.
	Since $0\le q\le1$, every probability measure $\mu$ satisfies
	\[
	\iint q(x-y)\,\mu(dx)\mu(dy)
	\le \sup_x\mu(B(x,R))+\frac{\delta}{2}.
	\]
	Each $B(x,R)$ lies in some $B(z,r)$ with $z\in\Z^d$ and
	$r=R+\sqrt d$. Hence, on
	$\{I_\sigma^q\ge\delta,\ Z_\sigma^\nu>\eta u\}$, we have
	\[
	\sup_{z\in\Z^d}M_\sigma^\nu(B(z,r))
	=Z_\sigma^\nu\sup_{z\in\Z^d}\mu_\sigma(B(z,r))
	\ge \frac{\eta\delta u}{2}.
	\]
Using
	\eqref{eq:sharpness-exact-hit}, we obtain
	\[
	\Pp\biggl(\sup_{t\ge0,\,z\in\Z^d} M_t^\nu(B(z,r))
	\ge \frac{\eta\delta u}{2}\biggr)
	\ge \frac12\,\Pp(\tau_u<\infty)
	= \frac{1}{2u}.
	\]
	Setting $u=2v/(\eta\delta)$ gives \eqref{eq:1110a}
	with $c=\eta\delta/4$.

\end{proof}

The following lemma gives a lower bound on the probability of observing a large total mass at a deterministic time.

\begin{lemma}\label{lem:logtime}
Under the assumptions of Lemma~\ref{lem:excursion},
for every $\varepsilon>0$ there exists $C_\varepsilon<\infty$ such that, for every
sufficiently large $v$, there is a deterministic integer
$s \in [1, C_\varepsilon\log v]$ for which
\[
 \Pp(Z_s\ge v)\ge v^{-1-\varepsilon}.
\]
\end{lemma}

\begin{proof}
	Fix $\varepsilon>0$ and take $r,c$ from Lemma~\ref{lem:excursion}.
	Choose $b>1$ large enough that Lemma~\ref{lem:excursion} applies at level $2b$ and
	\[
	b^{-1-\varepsilon/2}\le p:=\frac{c}{4b}<1.
	\]
	Let $\mathcal P_r$ denote the set of probability measures on $\R^d$
	supported on $\overline{B(0,r)}$. For $D_m:=2^{-m}\mathbb N\cap(0,m]$, set
	\[
	F_m(\nu):=\Pp\biggl(
	\sup_{t\in D_m,\,z\in\Z^d}M_t^\nu(B(z,r))>b
	\biggr).
	\]
	Since $t\mapsto M_t^\nu(B(z,r))$ is lower semicontinuous and the
	dyadic times are dense, Lemma~\ref{lem:excursion} gives
	\[
	\lim_{m\to\infty}F_m(\nu)
	\ge \Pp\biggl(\sup_{t\ge0,\,z\in\Z^d}
	M_t^\nu(B(z,r))\ge2b\biggr)
	\ge 2p.
	\]
	For each $t\in D_m$, use the continuous version of
    $\mathcal K_{0,t}^\kappa(x,y)$ from the proof of Lemma~\ref{lem:derivatives}.
    For every $z\in\Z^d$, the function
    $x\mapsto\int_{B(z,r)}\mathcal K_{0,t}^\kappa(x,y)\,dy$ is then
    continuous on $\overline{B(0,r)}$, by uniform continuity of the kernel
    on $\overline{B(0,r)}\times\overline{B(z,r)}$.
    On a common event of probability one, this holds for all the grid
    times and centers. Thus $\nu\mapsto M_t^\nu(B(z,r))$ is continuous
    on $\mathcal P_r$ for every such pair. The strict superlevel event
    defining $F_m$ is open as a set of initial measures for each such
    environment, so Fatou's lemma gives lower semicontinuity of $F_m$. The open sets
	$\{F_m>p\}$ increase and cover $\mathcal P_r$. By compactness,
	we may fix $m$ such that
	\[
	\inf_{\nu\in\mathcal P_r}F_m(\nu)\ge p.
	\]
	
	For $k\ge1$, define
	\[
	p_k:=\inf_{\nu\in\mathcal P_r}
	\Pp\biggl(
	\sup_{\substack{t\in2^{-m}\mathbb N\cap(0,km]\\z\in\Z^d}}
	M_t^\nu(B(z,r))>b^k
	\biggr).
	\]
	
Set
	\[
	\tau:=\inf\bigl\{t\in D_m:
	\sup_{z\in\Z^d}M_t^\nu(B(z,r))>b\bigr\},
	\qquad \inf\varnothing:=\infty.
	\]
	Fix an enumeration $(z_j)_{j\ge1}$ of $\Z^d$. On $\{\tau<\infty\}$, let
	\[
	J:=\min\{j\ge1:M_\tau^\nu(B(z_j,r))>b\},
	\]
	and define the recentered probability measure by
	\[
	\widehat\nu(A):=
	\frac{M_\tau^\nu\bigl((A+z_J)\cap B(z_J,r)\bigr)}
	{M_\tau^\nu(B(z_J,r))},
	\qquad A\in\mathcal B(\R^d).
	\]
	Since $\widehat\nu(B(0,r))=1$, we have $\widehat\nu\in\mathcal P_r$.

	On $\{\tau<\infty\}$, let $\widehat M$ start from $\widehat\nu$
	and be driven by the shifted noise
	$W_{\tau+s}(z_J+x)-W_\tau(z_J+x)$. For each $t\in D_m$ and $j\ge1$, the event $\{\tau=t,J=j\}$ belongs to $\mathcal F_t$. On this event the future noise is shifted by the fixed vector $z_j$. Its conditional law given $\mathcal F_t$ is the original noise law, by independent time increments and spatial stationarity. Since $\widehat\nu\in\mathcal P_r$ is known at $\tau$, the restart property and the definition of $p_k$ give
	\[
	\Pp\biggl(
	\sup_{\substack{s\in 2^{-m}\mathbb N\cap(0,km]\\ z\in\Z^d}}
	\widehat M_s(B(z,r))>b^k
	\,\biggm|\,\mathcal F_\tau
	\biggr)\ge p_k.
	\]
	By positivity of the kernel and the definition of $\widehat\nu$,
	\begin{align*}
		M_{\tau+s}^\nu(B(z_J+z,r))
		&\ge \int_{B(z_J,r)} M_\tau^\nu(dx)
		\int_{B(z_J+z,r)}\mathcal K_{\tau,\tau+s}^\kappa(x,y)\,dy\\
		&= M_\tau^\nu(B(z_J,r))\,\widehat M_s(B(z,r)).
	\end{align*}
Suppose that
	$\widehat M_s(B(z,r))>b^k$ for some $s\in 2^{-m}\mathbb N\cap(0,km]$ and $z\in\Z^d$.
	For this pair $(s,z)$,
	$
	M_{\tau+s}^\nu(B(z_J+z,r))
	> b^{k+1}.
	$
	 Hence, for every
	$\nu\in\mathcal P_r$,
	\[
	\Pp\biggl(
	\sup_{\substack{t\in 2^{-m}\mathbb N\cap(0,(k+1)m]\\ z\in\Z^d}}
	M_t^\nu(B(z,r))>b^{k+1}
	\biggr)
	\ge p_k\,\Pp(\tau<\infty)
	=p_k F_m(\nu).
	\]
	Taking the infimum over $\nu\in\mathcal P_r$ yields
	\[
	p_{k+1}\ge p_k\inf_{\nu\in\mathcal P_r}F_m(\nu)\ge p\,p_k.
	\]

	For the process started from $\delta_0$, the total mass dominates the mass in
	each ball. Summing over the $2^mkm$ grid times gives
	\[
	p^k
	\le \sum_{t\in2^{-m}\mathbb N\cap(0,km]}\Pp(Z_t>b^k)
	\le 2^mkm\max_{t\in2^{-m}\mathbb N\cap(0,km]}\Pp(Z_t>b^k).
	\]
	Choose a maximizing $t$ and set $s=\lceil t\rceil$. Since $s-t<1$,
	the martingale property and \eqref{eq:sharpness-replica-finite} give
	\[
	\E[Z_s\mid\mathcal F_t]=Z_t,\qquad
	\E[Z_s^2\mid\mathcal F_t]\le e^\kappa Z_t^2.
	\]
	The conditional Paley--Zygmund inequality therefore yields
	\[
	\Pp(Z_s\ge b^k/2)
	\ge \frac{e^{-\kappa}}4\,\Pp(Z_t>b^k)
	\ge \frac{e^{-\kappa}p^k}{4\cdot2^mkm}.
	\]
	
	Take $k=\lceil\log_b(2v)\rceil$. Then $b^k/2\ge v$,
	$1\le s\le km\le C_\varepsilon\log v$, and
	$p^k\ge p(2v)^{-1-\varepsilon/2}$. Consequently,
	\[
	\Pp(Z_s\ge v)
	\ge \frac{e^{-\kappa}p(2v)^{-1-\varepsilon/2}}{4\cdot2^mkm}
	\ge v^{-1-\varepsilon}
	\]
	for all sufficiently large $v$.
\end{proof}

\begin{lemma}\label{lem:poly}
Assume $\kappa>b_*$ and strong disorder. Then, for every $\theta\in(0,1)$ and
every $K>0$, all sufficiently large integers $n$ satisfy
\[
 \E Z_n^\theta\le2n^{-K}.
\]
\end{lemma}

\begin{proof}
	Fix $\theta\in(0,1)$ and $K>0$. Choose
	\[
	L>d+1+\frac{K}{1-\theta},
	\qquad u=n^L,\qquad \varepsilon=\frac1{2L}.
	\]
	By Lemma~\ref{lem:logtime}, there is an integer $1\le s\le C\log n$
	such that $\Pp(Z_s\ge u)\ge u^{-1-\varepsilon}$. Hence
	\[
	p:=\E[Z_s\mathbf1_{\{Z_s\ge u\}}]
	\ge u\,\Pp(Z_s\ge u)\ge n^{-1/2}.
	\]
	
	Following the block test of Junk and Lacoin~\cite[Lemma~3.2 and Proposition~4.1]{JL}, set
	$m=\lfloor n/s\rfloor$ and $Q_n=[-n,n]^d$, and let $U_j(x)$ denote
	the total mass started from $\delta_x$ using the noise on
	$[js,(j+1)s]$. Set $A_0=\Omega$ and
	\[
	A_{j+1}:=A_j\cap\bigl\{\sup_{x\in Q_n}U_j(x)<u\bigr\},
	\qquad 0\le j<m.
	\]
	Cover $Q_n$ by at most $C n^d$ translates of $[0,1]^d$. By
    \eqref{eq:sharpness-derivative-sup} with exponent one and spatial
    stationarity, $\E\sup_{x\in Q_n}Z_s(x)\le Cn^d$. Hence the union
    bound and Markov's inequality give
	\[
	\Pp(A_m^c)
	\le \frac{m}{u}\,\E\sup_{x\in Q_n}Z_s(x)
	\le C n^{1+d-L}.
	\]
	
	We next estimate $a_j:=\E[Z_{js}\mathbf1_{A_j}]$.
	For $x\in Q_n$, spatial stationarity gives
	\[
	\E\bigl[U_j(x)\mathbf1_{\{\sup_{y\in Q_n}U_j(y)<u\}}\bigr]
	\le \E[Z_s\mathbf1_{\{Z_s<u\}}]=1-p.
	\]
	For $x\notin Q_n$, the same expectation is at most $\E[U_j(x)]=1$.
	Since $A_j\in\mathcal F_{js}$ and the next block is independent of
	$\mathcal F_{js}$, the flow identity
	\[
	Z_{(j+1)s}=\int U_j(x)\,M_{js}^{\delta_0}(dx)
	\]
	gives
	\[
	a_{j+1}
	\le (1-p)a_j+p\,\E[M_{js}^{\delta_0}(Q_n^c)].
	\]
	By the Feynman--Kac formula and the Gaussian tail bound,
	\[
	\E[M_{js}^{\delta_0}(Q_n^c)]
	=P_0^B(B_{js}\notin Q_n)\le 2d e^{-n/2},
	\qquad js\le n.
	\]
	Iterating from $a_0=1$ and using $m\ge cn/\log n$, we obtain
	\[
	a_m\le (1-p)^m+2d e^{-n/2}
	\le \exp\{-c\sqrt n/\log n\}+2d e^{-n/2}.
	\]
	
	Finally, $A_m\in\mathcal F_{ms}$, so the martingale property gives
	$\E[Z_n\mathbf1_{A_m}]=a_m$. Applying H\"older's inequality on
	$A_m^c$ and $A_m$ and using $\E Z_n=1$ yields
	\[
	\E[Z_n^\theta]
	\le \Pp(A_m^c)^{1-\theta}+a_m^\theta
	\le 2n^{-K}
	\]
	for all sufficiently large $n$, since $(L-d-1)(1-\theta)>K$.
\end{proof}

\begin{lemma}\label{lem:finitevolume}
	Fix $\kappa>0$ and $\theta\in(0,1)$. There exist $K_0,n_0<\infty$
	such that if $\E Z_n^\theta\le 2n^{-K_0}$ for some integer
	$n\ge n_0$, then
	\[
	\E Z_t^\theta\le Ce^{-ct},\qquad t\ge0,
	\]
	for some $C,c>0$. The same conclusion holds whenever
	$S_n(\kappa)<1$ for some integer $n\ge1$, where $S_n(\kappa)$
	is defined in \eqref{eq:sharpness-cell-sum}.
\end{lemma}

\begin{proof}
	Following \cite[Proposition~3.1]{JL}, let
	$\mathcal C_z=z+[0,1)^d$ and define
	\begin{equation}\label{eq:sharpness-cell-sum}
		A_n^\kappa(z,z')
		:=\sup_{x\in\overline{\mathcal C_z}}
		\int_{\mathcal C_{z'}}\mathcal K_{0,n}^\kappa(x,y)\,dy,
		\qquad
		S_n(\kappa):=\sum_{z\in\Z^d}\E[A_n^\kappa(0,z)^\theta].
	\end{equation}
	Write $A_n^{(j)}(z,z')$ for the same quantity on $[jn,(j+1)n]$.
	By kernel composition,
	\[
	M_{(j+1)n}^{\delta_0}(\mathcal C_{z'})
	\le \sum_{z\in\Z^d}
	M_{jn}^{\delta_0}(\mathcal C_z)\,A_n^{(j)}(z,z').
	\]
	Set $T_j:=\sum_z\E[M_{jn}^{\delta_0}(\mathcal C_z)^\theta]$.
	By subadditivity, independence of the next time block, and spatial
	stationarity,
	\[
	T_{j+1}
	\le \sum_z\E[M_{jn}^{\delta_0}(\mathcal C_z)^\theta]
	\sum_{z'}\E[A_n^\kappa(z,z')^\theta]
	=T_j S_n(\kappa).
	\]
	Since $T_0=1$, another application of subadditivity gives
	\begin{equation}\label{eq:sharpness-finite-volume-product}
		\E Z_{jn}^\theta\le T_j\le S_n(\kappa)^j.
	\end{equation}
	If $S_n(\kappa)<1$, the supermartingale property of $Z^\theta$ gives
	\[
	\E Z_t^\theta
	\le S_n(\kappa)^{\lfloor t/n\rfloor}
	\le S_n(\kappa)^{-1}
	\exp\bigl\{t\log S_n(\kappa)/n\bigr\}.
	\]
	
	To verify the moment criterion, write $a_z=A_n^\kappa(0,z)$ and
    $Q=[0,1]^d$. Since $a_z\le\sup_{x\in Q}Z_n(x)$, the unit-cube
    comparison \eqref{eq:sharpness-derivative-sup} gives
    $\E a_z^\theta\le C_\theta\E Z_n^\theta$. On the other hand,
    Jensen's inequality and \eqref{eq:sharpness-kernel-tail} give
    \[
      \E a_z^\theta\le(\E a_z)^\theta,\qquad
      \E a_z\le C\int_{\mathcal C_z}p_{n+1}(y)\,dy
                 \le C n^{-d/2}e^{-c|z|^2/n}.
    \]
    Taking the geometric mean of the two bounds on $\E a_z^\theta$
    and summing over $z$ yields
    \[
      S_n(\kappa)^2
      \le C_\theta\E Z_n^\theta
           \left[\sum_{z\in\Z^d}(\E a_z)^{\theta/2}\right]^2
      \le C_\theta n^{d(1-\theta/2)}\E Z_n^\theta .
    \]
    Choose $K_0>d(1-\theta/2)$. If $\E Z_n^\theta\le2n^{-K_0}$,
    then $S_n(\kappa)<1$ for all sufficiently large such $n$.
    The conclusion follows from \eqref{eq:sharpness-finite-volume-product}.
\end{proof}

\begin{proof}[Proof of Theorem~\ref{thm:sharp}]
	Assume that $\kappa>b_*$ and $Z_t\to0$ almost surely. Fix $\theta\in(0,1)$.
	Let $K_0$ and $n_0$ be as in Lemma~\ref{lem:finitevolume}.
	By Lemma~\ref{lem:poly}, we may choose an integer $n\ge n_0$
	such that $\E Z_n^\theta\le 2n^{-K_0}$. Lemma~\ref{lem:finitevolume}
	then gives constants $C,c>0$ such that
	\[
	\E Z_t^\theta\le Ce^{-ct},\qquad t\ge0.
	\]
	
	By the Feynman--Kac representation and Jensen's inequality,
	\[
	\log Z_t\ge \beta E_0^B[H_t(B)]-\frac{\kappa t}{2}.
	\]
	The right-hand side is integrable, and $\log^+ Z_t\le Z_t$.
	Thus $\log Z_t$ is integrable. Applying Jensen's inequality to
	$Z_t^\theta$ now yields
	\[
	\E\log Z_t
	\le \frac1\theta\log\E Z_t^\theta
	\le \frac{\log C-ct}{\theta}.
	\]
	Dividing by $t$ and letting $t\to\infty$, we obtain
	\[
	\limsup_{t\to\infty}\frac1t\E\log Z_t
	\le -\frac c\theta<0.
	\]
	
	If $\kappa>\kappa_c$, then $\kappa>b_*$ by Lemma~\ref{lem:spectral},
	and Propositions~\ref{prop:critical-parameter} and~\ref{prop:mean-dichotomy}
	give $Z_t\to0$ almost surely. Hence both conclusions hold.
\end{proof}

\section{Exponential fractional moments imply local extinction}
\label{sec:eventual-extinction}

Fix $\kappa$ and write $Z_t(x)=Z_t^\kappa(x)$. We prove local
extinction under the assumption that, for some $\theta\in(0,1)$
and constants $b,C>0$,
\begin{equation}\label{fracdecay}
	\E[Z_t(0)^\theta]\le Ce^{-bt},\qquad t\ge1.
\end{equation}
By Theorem~\ref{thm:sharp}, this bound holds whenever $\kappa>\kappa_c$.
We use the decomposition of the initial measure into its restrictions
to a ball and its complement, as in Mytnik and Xiong~\cite[Section~4]{MX}
and Chen, Ren and Zhao~\cite[proof of Theorem~1.5]{CRZ}. Taking the
radius proportional to time gives exponential bounds both for survival
of the inner mass and for occupation of a fixed compact region by the
outer mass.

We first extend the fractional-moment estimate \eqref{fracdecay} to
Lebesgue initial measure restricted to a ball whose radius grows
linearly with time.

\begin{lemma}\label{lem:expanding-ball}
For each $A>0$, there are $b_1,C_A>0$ such that
\begin{equation}\label{ballmoment}
 \E\left[\left(\int_{B(0,At)}Z_t(x)\,dx\right)^\theta\right]
       \le C_A t^d e^{-b_1t},\qquad t\ge1.
\end{equation}
\end{lemma}
\begin{proof}
Let $Q=[0,1]^d$ and $J_{A,t}=\{j\in\mathbb Z^d:(j+Q)\cap B(0,At)\ne\varnothing\}$.
For $t\ge1$, $|J_{A,t}|\le C_A t^d$. Since each cube has volume one,
subadditivity of $u\mapsto u^\theta$ and spatial stationarity give
\begin{equation}\label{eq:ball-cube-cover}
 \E\left[\left(\int_{B(0,At)}Z_t(x)\,dx\right)^\theta\right]
 \le \sum_{j\in J_{A,t}}\E\left[\left(\sup_{x\in j+Q}Z_t(x)\right)^\theta\right]
 = |J_{A,t}|\,\E\left[\left(\sup_{x\in Q}Z_t(x)\right)^\theta\right].
\end{equation}
By \eqref{eq:sharpness-derivative-sup} and \eqref{fracdecay}, the last
expectation is at most $C_\theta\E[Z_t(0)^\theta]\le C_\theta C e^{-bt}$.
This proves \eqref{ballmoment} with $b_1=b$.
\end{proof}

We next use Lemma~\ref{lem:expanding-ball} to bound the probability
that the superprocess started from $m|_{B(0,AT)}$ survives until time $2T$.

\begin{lemma}\label{lem:inner-survival}
For each fixed $A>0$, there are $c,C_A>0$ such that
\[
 \Pp_{m|_{B(0,AT)}}(X_{2T}\ne0)\le C_Ae^{-cT},\qquad T\ge1.
\]
\end{lemma}
\begin{proof}
	
	Fix $\eta\in\Mf$ and $T>0$. For $k>0$, let $u^k$ solve the nonlinear
	equation in \eqref{eq:dual-equations} with constant initial value $k$.
	By spatial stationarity, the first and second moments of $u_t^k(x)$
	are independent of $x$. Write $m_k(t):=\E u_t^k(x)$. By comparison,
	$u_s^k\le k v_s^1$. The two-bridge estimate preceding
	\eqref{eq:integrated-dual}, applied with initial value $1$, gives
	$\E[(v_s^1(x))^2]\le e^{\kappa s}$. In particular,
	\[
	   \sup_{0\le s\le T,\,x\in\R^d}\E[(u_s^k(x))^2]\le k^2 e^{\kappa T}.
	\]
	Since $q\le1$ and $p_{t-s}$ has integral one, Cauchy--Schwarz bounds
	the expected bracket of the mild stochastic convolution by
	$\kappa k^2\int_0^t e^{\kappa s}\,ds$. It therefore has mean zero.
	The heat semigroup preserves constants, so taking expectations gives
	\[
	m_k(t)=k-\frac12\int_0^t \E[(u_s^k(x))^2]\,ds.
	\]
	Hence, by Jensen's inequality,
	$m_k'(t)\le -m_k(t)^2/2$. Since $m_k(0)=k$ and $(1/m_k)' \ge 1/2$ almost everywhere on
	$\{m_k>0\}$, integration gives
	\[
	m_k(t)\le \frac{1}{k^{-1}+t/2}=\frac{k}{1+kt/2}.
	\]
	The bound also holds when $m_k(t)=0$.
	
The annealed form of \eqref{eq:quenched-reversed-laplace} with $f\equiv k$ gives
\[
\E_\eta e^{-kX_T(1)}=\E e^{-\eta(u_T^k)}.
\]
By spatial stationarity and Tonelli's theorem, we have
$$\E[\eta(u_T^k)]=\int \E[u_T^k(x)]\,\eta(dx)=\eta(1)m_k(T).$$
Jensen's inequality and the upper bound on $m_k(T)$ yield
\[
\E e^{-\eta(u_T^k)}
\ge e^{-\eta(1)m_k(T)}
\ge \exp\biggl\{-\frac{k\eta(1)}{1+kT/2}\biggr\}.
\]
	Since $e^{-kX_T(1)}\downarrow\mathbf1_{\{X_T=0\}}$ as $k\to\infty$,
	bounded convergence yields
	\begin{equation}\label{eq:annealed-zero-bound}
		\Pp_\eta(X_T=0)
		=\lim_{k\to\infty}\E_\eta e^{-kX_T(1)}
		\ge \lim_{k\to\infty}\exp\biggl\{-\frac{k\eta(1)}{1+kT/2}\biggr\}
		=e^{-2\eta(1)/T}.
	\end{equation}

By \eqref{eq:first-moment}, $\E_\eta X_T(1)=\eta(1)<\infty$.
Applying the Markov property at $T$ and then \eqref{eq:annealed-zero-bound}
to the finite measure $X_T$ gives
\begin{align*}
 \Pp_\eta(X_{2T}=0)
 &=\E_\eta\big[\Pp_{X_T}(X_T=0)\big]\\*
 &\ge\E_\eta e^{-2X_T(1)/T}
 =\E e^{-\eta(u_T^{2/T})}
 \ge\E e^{-(2/T)\eta(v_T^1)}.
\end{align*}
Thus
\begin{equation}\label{eq:two-stage-survival}
 \Pp_\eta(X_{2T}\ne0)
       \le\E[1-e^{-(2/T)\eta(v_T^1)}].
\end{equation}

For the linear kernel, set
$W^{\leftarrow T}(s,x)=W(T,x)-W(T-s,x)$ for $0\le s\le T$.
In the Brownian bridge representation from the proof of
Lemma~\ref{lem:derivatives}, reversing the bridge and the order of
the environmental increments gives
\[
\mathcal K_{0,T}^{W}(x,y)
=\mathcal K_{0,T}^{W^{\leftarrow T}}(y,x).
\]
Integrating over $x$ and using the kernel formulas in
Section~\ref{sec:bridge}, we obtain
\[
v_T^1(y)
=\int_{\R^d}\mathcal K_{0,T}^{W^{\leftarrow T}}(y,x)\,dx
=Z_T^{W^{\leftarrow T}}(y).
\]
Since $W^{\leftarrow T}$ has the same law as $W$, it follows that
\[
\int_{B(0,AT)}v_T^1(y)\,dy
\ \stackrel{\rm law}{=}\
\int_{B(0,AT)}Z_T(y)\,dy.
\]

Taking $\eta=m|_{B(0,AT)}$ in \eqref{eq:two-stage-survival} and using
$1-e^{-x}\le x^\theta$ and \eqref{ballmoment}, we obtain
\begin{align*}
 \Pp_{m|_{B(0,AT)}}(X_{2T}\ne0)
 &\le\E\left[1-\exp\left\{-\frac2T
                  \int_{B(0,AT)}Z_T(x)\,dx\right\}\right]\\
 &\le(2/T)^\theta
          \E\left(\int_{B(0,AT)}Z_T(x)\,dx\right)^\theta\\
 &\le C_A T^{d-\theta}e^{-b_1T}
 \le C_A e^{-b_1T/2}.
\end{align*}
\end{proof}

We next estimate the probability that mass initially outside \(B(0,An)\) reaches the fixed region \(\{f>0\}\) by time \(n+2\).

\begin{lemma}\label{lem:outer-occupation}
Let $f\in C_c^{\infty,+}$ and $A>0$. There are $c,C>0$ such that,
for every sufficiently large integer $n$,
\begin{equation}\label{eq:outer-occupation}
 \Pp_{m|_{B(0,An)^c}}
       \left(\int_0^{n+2}X_s(f)\,ds>0\right)\le Ce^{-cn}.
\end{equation}
\end{lemma}

\begin{proof}
	Choose $r>0$ so that $\operatorname{supp}f\subset B(0,r)$, and write
	\[
	Y_T(f):=\int_0^T X_s(f)\,ds,
	\qquad
	\mu_n:=m|_{B(0,An)^c}.
	\]
	For $\ell>0$, let $v^\ell$ denote the nonnegative solution of
	\[
	\partial_t v^\ell
	=\tfrac12\Delta v^\ell-\tfrac12(v^\ell)^2+\ell f,
	\qquad v_0^\ell=0.
	\]
	Applying \cite[Theorem~2.18]{MX}, we obtain
	\[
	\E_{\mu_n}e^{-\ell Y_T(f)}
	\ge \exp\{-\mu_n(v_T^\ell)\}.
	\]
	Consequently,
	\begin{equation}\label{eq:occupation-probability-bound}
		\E_{\mu_n}\bigl[1-e^{-\ell Y_T(f)}\bigr]
		\le 1-e^{-\mu_n(v_T^\ell)}
		\le \mu_n(v_T^\ell)
		= \int_{|x|\ge An}v_T^\ell(x)\,dx .
	\end{equation}
	We next estimate the last integral uniformly in $\ell$.
	
	First, we show that $v_t^\ell(x)\le 6$ on the sphere $|x|=r+1$,
	uniformly in $t$ and $\ell$. For $|x|>r$, set
	\[
	w(x):=\frac{6}{(|x|-r)^2}.
	\]
	A direct calculation gives
	\[
	\Delta w(x)
	=\frac{36}{(|x|-r)^4}
	-\frac{12(d-1)}{|x|(|x|-r)^3}
	\le w(x)^2.
	\]
	Thus $w$ is a stationary supersolution of the equation for
	$v^\ell$ on $\{|x|>r\}$, where $f=0$.
	
Note that
	\begin{equation}\label{lm}
	0\le v_t^\ell(x)
	\le \ell\int_0^t(p_s*f)(x)\,ds
	\le \ell T\|f\|_\infty,
	\qquad 0\le t\le T.
	\end{equation}
	Since $f$ has compact support, the middle expression tends to zero
	as $|x|\to\infty$, uniformly for $0\le t\le T$. 
	
	On the inner boundary $|x|=r+\varepsilon$, we have
	$w(x)=6/\varepsilon^2$.  By \eqref{lm}, we may choose $\varepsilon\in(0,1)$ small
	enough that
	\[
	v_t^\ell(x)\le\ell T\|f\|_\infty
	\le \frac{6}{\varepsilon^2}=w(x)
	\qquad\text{for } |x|=r+\varepsilon,\quad 0\le t\le T.
	\]
	For the outer boundary, fix $\delta>0$. By the uniform decay
	established above, we may choose $R>r+1$ so large that
	$v_t^\ell(x)\le\delta$ whenever $|x|=R$ and $0\le t\le T$.

	The function $w+\delta$ is also a supersolution. It dominates
	$v^\ell$ at time zero and on both boundary spheres of the annulus
	$r+\varepsilon<|x|<R$. Parabolic comparison therefore gives
	\[
	v_t^\ell(x)\le w(x)+\delta=6+\delta,
	\qquad |x|=r+1,\quad 0\le t\le T.
	\]
	Letting $\delta\downarrow0$ proves the asserted boundary bound.
	
	Next, put $R_0=r+1$.
	For each $\lambda>0$ and $|x|>R_0$, define
	\[
	h_\lambda(t,x)
	:=6\exp\bigl\{\tfrac12\lambda^2t-\lambda(|x|-R_0)\bigr\}.
	\]
	A direct calculation gives
	\[
	(\partial_t-\tfrac12\Delta)h_\lambda(t,x)
	=\frac{(d-1)\lambda}{2|x|}\,h_\lambda(t,x)\ge0.
	\]
On $|x|=R_0$ we have $v_t^\ell\le6\le h_\lambda(t,x)$. For fixed
$T,\ell$ and $\delta>0$, \eqref{lm} gives $v_t^\ell\le\delta$ at every
sufficiently large outer radius, uniformly for $t\le T$. Compare the
heat subsolution $v^\ell$ with $h_\lambda+\delta$ on these annuli,
let the outer radius tend to infinity, and then let $\delta\downarrow0$.
This gives
	\[
	v_T^\ell(x)
	\le 6\exp\bigl\{\tfrac12\lambda^2T-\lambda(|x|-R_0)\bigr\},
	\qquad |x|>R_0.
	\]
Choosing $\lambda=(|x|-R_0)/T$ yields
	\begin{equation}\label{eq:occupation-hit-bound}
		v_T^\ell(x)
		\le 6\exp\biggl\{-\frac{(|x|-R_0)^2}{2T}\biggr\},
		\qquad |x|>R_0,
	\end{equation}
	uniformly in $\ell>0$.
	
	Finally, take $T=n+2$ and choose $n$ large enough that $An\ge2R_0$.
	For $|x|\ge An$, we have $|x|-R_0\ge |x|/2$ and $n+2\le3n$.
	It follows from \eqref{eq:occupation-hit-bound} that
	\[
	\begin{aligned}
		\int_{|x|\ge An}v_{n+2}^\ell(x)\,dx
		&\le 6\int_{|x|\ge An}e^{-|x|^2/(24n)}\,dx\\
		&\le 6e^{-A^2n/48}
		\int_{\R^d}e^{-|x|^2/(48n)}\,dx\\
		&=6(48\pi n)^{d/2}e^{-A^2n/48}
		\le C e^{-A^2n/96}.
	\end{aligned}
	\]
	This bound is independent of $\ell$. Since $Y_{n+2}(f)\ge0$,
	\[
	1-e^{-\ell Y_{n+2}(f)}
	\uparrow \mathbf{1}_{\{Y_{n+2}(f)>0\}}
	\qquad\text{as }\ell\to\infty.
	\]
	Applying monotone convergence in
	\eqref{eq:occupation-probability-bound} now gives
	\[
	\Pp_{\mu_n}\bigl(Y_{n+2}(f)>0\bigr)
	\le C e^{-A^2n/96},
	\]
	which proves \eqref{eq:outer-occupation}.
\end{proof}

Lemmas~\ref{lem:inner-survival} and~\ref{lem:outer-occupation}
together imply that the mass of any fixed compact set is zero
for all sufficiently large times.

\begin{proposition}\label{prop:eventual-extinction}
Assume \eqref{fracdecay} with exponent $b>0$ and set
$c=\min\{b/4,1/96\}$. For every compact $K\subset\R^d$, there exist
$C_K<\infty$ and an almost surely finite random time $T_K$ such that
almost surely $X_s(K)=0$ for all $s\ge T_K$ and
\[
 \Pp_m(L_K>t)\le\Pp_m(T_K>t)\le C_K e^{-ct},\qquad t\ge0,
\]
where $L_K=\sup(\{s\ge0:X_s(K)>0\}\cup\{0\})$.
\end{proposition}

\begin{proof}
	Fix a compact set $K\subset\mathbb R^d$ and choose
	$f\in C_c^{\infty,+}$ such that $f\ge1$ on $K$. Then
	$X_t(K)\le X_t(f)$ for every $t\ge0$. For each integer $n\ge1$, set
	\[
	E_n:=\biggl\{\sup_{n\le t\le n+1}X_t(f)>0\biggr\}.
	\]
	We first bound $\Pp_m(E_n)$ for large $n$.
	
	Fix $n$ and decompose the initial measure as
	\[
	m=\mu_n+\nu_n,\qquad
	\mu_n:=m|_{B(0,n)},\qquad
	\nu_n:=m|_{B(0,n)^c}.
	\]
	By the conditional branching property recalled in
	Section~\ref{sec:main}, the process under $\Pp_m$ has the same law
	as the sum of two processes $X^{\mathrm{in}}$ and $X^{\mathrm{out}}$
	started from $\mu_n$ and $\nu_n$, respectively, in the same
	environment. Since zero is absorbing, $X_n^{\mathrm{in}}=0$ implies
	$X_t^{\mathrm{in}}=0$ for all $t\ge n$. Moreover, if
	$X_t^{\mathrm{out}}(f)>0$ at some $t\in[n,n+1]$, then the continuity
	and nonnegativity of $s\mapsto X_s^{\mathrm{out}}(f)$ imply that
	\[
	\int_0^{n+2}X_s^{\mathrm{out}}(f)\,ds>0.
	\]
	
	Let $\mathbf Q_n$ denote the joint law of
	$X^{\mathrm{in}}$ and $X^{\mathrm{out}}$ in the common environment.
	Since their sum has law $\Pp_m$, the preceding event inclusion
	and the union bound give
	\[
	\begin{aligned}
		\Pp_m(E_n)
		&= \mathbf Q_n\biggl(
		\sup_{n\le t\le n+1}
		(X_t^{\mathrm{in}}+X_t^{\mathrm{out}})(f)>0
		\biggr)\\
		&\le \mathbf Q_n(X_n^{\mathrm{in}}\ne0)
		+ \mathbf Q_n\biggl(\int_0^{n+2}X_s^{\mathrm{out}}(f)\,ds>0\biggr)\\
		&= \Pp_{\mu_n}(X_n\ne0)
		+ \Pp_{\nu_n}\biggl(\int_0^{n+2}X_s(f)\,ds>0\biggr).
	\end{aligned}
	\]

Take $T=n/2$ and $A=2$ in the proof of Lemma~\ref{lem:inner-survival}.
Since $b_1=b$ in Lemma~\ref{lem:expanding-ball}, the first probability
is at most $C n^{d-\theta}e^{-bn/2}\le C'e^{-bn/4}$ for $n\ge2$. The proof of
Lemma~\ref{lem:outer-occupation} with $A=1$ bounds the second by
$C e^{-n/96}$ for all sufficiently large $n$, with the starting index
depending on the support of $f$. Thus $\Pp_m(E_n)\le C_K e^{-cn}$
for all $n\ge1$, after increasing $C_K$ to cover the remaining
integers. Set
\[
 T_K:=1+\sup\bigl(\{n\ge1:E_n\text{ occurs}\}\cup\{0\}\bigr).
\]
For $t\ge1$, the union bound gives
\[
 \Pp_m(T_K>t)
 \le \sum_{n\ge\lfloor t\rfloor}\Pp_m(E_n)
 \le \frac{C_K e^{-c\lfloor t\rfloor}}{1-e^{-c}} .
\]
This yields the asserted exponential bound for $t\ge0$ after an adjustment of the constant. In particular, $T_K<\infty$ almost surely. On this event $X_s(f)=0$ on every interval $[n,n+1]$ with $n\ge T_K$, so $X_s(K)=0$ for all $s\ge T_K$ and $L_K\le T_K$.

For completeness, choose $g_j\in C_c^+$ decreasing pointwise to $\mathbf1_K$. For each compact time interval $I$, continuity of $s\mapsto X_s(g_j)$ and compactness of $I$ give
\[
 \sup_{s\in I}X_s(K)=\inf_j\max_{s\in I}X_s(g_j).
\]
Indeed, if $a$ is smaller than the right-hand side, the nested closed sets $\{s\in I:X_s(g_j)\ge a\}$ have nonempty intersection. Each maximum is measurable by continuity. Thus $\{L_K>t\}$ is measurable as the union of $\{\sup_{s\in I}X_s(K)>0\}$ over closed intervals $I\subset(t,\infty)$ with rational endpoints. The tail bound for $L_K$ follows.
\end{proof}

\section{Critical persistence}\label{sec:critical}
\subsection{Free energy at the critical parameter}
We first show that the free energy and the exponential rates of the fractional
moments vanish at $\kappa_c$. Together with Theorem~\ref{thm:sharp}, this
reduces uniform integrability at $\kappa_c$ to the strict inequality
$b_*<\kappa_c$. We establish this inequality below by proving weak disorder
on an interval extending to the right of the second-moment threshold $b_*$.

\begin{lemma}\label{lem:free-energy}
For every $\kappa\ge0$, the limit
\[
p(\kappa):=\lim_{T\to\infty}\frac1T\E\log Z_T^\kappa(0)
\]
exists and belongs to $[-\kappa/2,0]$. The function $p$ is
nonincreasing and $1/2$-Lipschitz. For each $\kappa\ge0$,
\[
\lim_{T\to\infty}\frac1T\log Z_T^\kappa(0)=p(\kappa)
\qquad\text{almost surely}.
\]
\end{lemma}
\begin{proof}
Write $a_T(\kappa)=\E\log Z_T^\kappa(0)$. With $Q=q$, $q(0)=1$ and inverse temperature $\beta=\sqrt\kappa$, the normalized partition function in \cite[Proposition~1.1]{Lacoin11} is precisely \eqref{eq:polymer-definition} at $x=0$. That proposition gives the deterministic limit of $a_T(\kappa)/T$ and the same almost sure limit of $T^{-1}\log Z_T^\kappa(0)$ for $\kappa>0$. At $\kappa=0$ both quantities vanish.

Fix $0\le\kappa\le\kappa'$ and put $\delta=\kappa'-\kappa$. Couple the environments as in Proposition~\ref{prop:critical-parameter}, using independent fields $W$ and $W'$ with covariance $q$. For fixed $W$, define the polymer expectation of an integrable path functional $F$ by
\[
 \langle F\rangle_{T,\kappa}
 =\frac{E_0^B\!\left[F(B)\exp\{\sqrt\kappa H_T(B)-\kappa T/2\}\right]}
        {Z_T^\kappa(0)}.
\]
Here $E_0^B$ averages over Brownian motion started at $0$. Put $H'_T(B)=\int_0^T W'(ds,B_s)$. Under this coupling,
\[
 R_T:=\frac{Z_T^{\kappa'}(0)}{Z_T^\kappa(0)}
 =\left\langle e^{\sqrt\delta H'_T(B)-\delta T/2}\right\rangle_{T,\kappa}.
\]
Write $\E'$ for expectation over $W'$ with $W$ fixed. For every fixed path, $H'_T(B)$ is centered Gaussian with variance $T$. Hence $\E'R_T=1$ and $\E'\langle |H'_T(B)|\rangle_{T,\kappa}=\sqrt{2T/\pi}$. Jensen's inequality under the path measure gives $\log R_T\ge\sqrt\delta\langle H'_T(B)\rangle_{T,\kappa}-\delta T/2$. This lower bound and $\log^+R_T\le R_T$ give integrability. Taking $\E'$ in the lower bound and applying Jensen's inequality under $\E'$ for the upper bound yields
\[
 -\frac{\delta T}{2}\le\E'\log R_T\le\log\E'R_T=0.
\]
Averaging over $W$, dividing by $T$ and letting $T\to\infty$, we obtain
\[
 -\tfrac12(\kappa'-\kappa)\le p(\kappa')-p(\kappa)\le0.
\]
Since $p(0)=0$, this also gives $-\kappa/2\le p(\kappa)\le0$.
\end{proof}

\begin{proposition}[Critical free energy]\label{prop:critical-free-energy}
The free energy vanishes at the critical parameter:
$p(\kappa_c)=0$. Moreover,
\[
\Pp\left(
\lim_{T\to\infty}\frac1T\log Z_T^{\kappa_c}(0)=0
\right)=1.
\]
For every $\theta\in(0,1)$, the fractional moment satisfies
\[
\lim_{T\to\infty}\frac1T
\log\E[(Z_T^{\kappa_c})^\theta]=0.
\]
\end{proposition}

\begin{proof}
	For $0\le\kappa<\kappa_c$,
	Propositions~\ref{prop:critical-parameter}
	and~\ref{prop:mean-dichotomy} give
	$
	\E Z_\infty^\kappa=c_\kappa=1,
$
so $
	\Pp(Z_\infty^\kappa>0)>0.
	$
By Lemma~\ref{lem:free-energy}, $p(\kappa)\le0$. If $p(\kappa)<0$, the same lemma gives $Z_T^\kappa \to 0$ almost surely, so $\Pp(Z_\infty^\kappa>0)=0$, a contradiction. Hence
	$
p(\kappa)=0
	$
	for
	$0\le\kappa<\kappa_c.$
Continuity of $p$ from Lemma~\ref{lem:free-energy} yields
	\[
	p(\kappa_c)=\lim_{\kappa\uparrow\kappa_c}p(\kappa)=0.
	\]
	The almost sure limit in the same lemma then gives
	\[
	\lim_{T\to\infty}\frac1T\log Z_T^{\kappa_c}(0)
	=p(\kappa_c)=0
	\qquad\text{almost surely}.
	\]
	
	For $\theta\in(0,1)$, Jensen's inequality and
	$\E Z_T^{\kappa_c}=1$ imply
	\[
	\frac{\theta}{T}\E\log Z_T^{\kappa_c}
	\le
	\frac1T\log\E[(Z_T^{\kappa_c})^\theta]
	\le
	\frac{\theta}{T}\log\E Z_T^{\kappa_c}
	=0.
	\]
	Since
	\[
	\lim_{T\to\infty}
	\frac{\theta}{T}\E\log Z_T^{\kappa_c}
	=\theta p(\kappa_c)=0,
	\]
	we conclude that
	\[
	\lim_{T\to\infty}\frac1T
	\log\E[(Z_T^{\kappa_c})^\theta]=0.
	\]
\end{proof}

\begin{theorem}\label{thm:gap-case}
	If $\kappa_c>b_*$, then the martingale
	$(Z_t^{\kappa_c}(0))_{t\ge0}$ is uniformly integrable.
\end{theorem}
\begin{proof}
	Suppose that $(Z_t^{\kappa_c}(0))_{t\ge0}$ is not uniformly
	integrable. Proposition~\ref{prop:mean-dichotomy} then gives
	\[
	Z_t^{\kappa_c}(0)\longrightarrow0
	\qquad\text{almost surely}.
	\]
	Since $\kappa_c>b_*$, Theorem~\ref{thm:sharp} implies
	\[
	p(\kappa_c)
	=\lim_{t\to\infty}\frac1t\E\log Z_t^{\kappa_c}(0)<0,
	\]
	contradicting Proposition~\ref{prop:critical-free-energy}.
	Hence $(Z_t^{\kappa_c}(0))_{t\ge0}$ is uniformly integrable.
\end{proof}

\subsection{The second-moment threshold and a scalar Gaussian projection}

Write $P_t=P_t^D$. Recall the strictly positive eigenfunction
$\psi\in C_0(\R^d)$ constructed in the proof of Lemma~\ref{lem:spectral}.
It satisfies
\begin{equation}\label{oo}
\psi\ge0,\qquad \psi\not\equiv0,\qquad
\psi=b_*T\psi=b_*G(q\psi).
\end{equation}
It follows that $(\Delta+b_*q)\psi=0$
in the sense of distributions.
The smoothness of $q$ and local elliptic regularity then give
$\psi\in C^\infty(\R^d)$. Consequently, the equation holds
pointwise:
\begin{equation}\label{r2}
\Delta\psi(x)+b_*q(x)\psi(x)=0,
\qquad x\in\R^d.
\end{equation}

\begin{proposition}[Exact second-moment threshold]\label{prop:exact-L2}
For every $\kappa\ge0$,
\[
\sup_{t\ge0}\E[(Z_t^\kappa(0))^2]<\infty
\quad\Longleftrightarrow\quad
\kappa<b_*.
\]
For every $\kappa\ge b_*$,
\[
\lim_{t\to\infty}\E[(Z_t^\kappa(0))^2]=\infty.
\]
\end{proposition}

\begin{proof}
	By the second-moment formula \eqref{eq:sharpness-replica-finite}
	with $\nu=\delta_0$, we have 
	\[
	\E[(Z_t^\kappa(0))^2]
	=E_0^D\!\left[
	\exp\left\{\kappa\int_0^t q(D_s)\,ds\right\}
	\right].
	\]
For $0\le\kappa<b_*$, the estimate in the proof of Lemma~\ref{lem:spectral} gives
	\[
	\sup_{t\ge0}\E[(Z_t^\kappa(0))^2]
	\le E_0^D\!\left[
	\exp\left\{\kappa\int_0^\infty q(D_s)\,ds\right\}
	\right]
	<\infty.
	\]
	
Conversely, suppose that $\kappa\ge b_*$. Note that
	\[
	\begin{aligned}
		\kappa^n(T^{n-1}Gq)(0)
		&\ge
		\frac{\kappa^n}{b_*\|\psi\|_\infty}
		(T^{n-1}\psi)(0)=\left(\frac{\kappa}{b_*}\right)^n
		\frac{\psi(0)}{\|\psi\|_\infty}
		\ge\frac{\psi(0)}{\|\psi\|_\infty}>0.
	\end{aligned}
	\]
Therefore,
	\[
	\lim_{t\to\infty}\E[(Z_t^\kappa(0))^2]
	=
	E_0^D\!\left[
	\exp\left\{\kappa\int_0^\infty q(D_s)\,ds\right\}
	\right]
	=\infty.
	\]
	Consequently, bounded second moments require $\kappa<b_*$,
	which proves the converse. 
	
\end{proof}

\begin{proposition}[Obstruction at the spectral endpoint]
\label{prop:auxiliary-endpoint-obstruction}
Assume that $\kappa=b_*$. Then there do not exist a function
$h\in C_0(\R^d)\cap C_b^2(\R^d)$ and a constant $a>0$ such that
\[
\Delta h(x)+\kappa q(x)h(x)\ge a\,q(x),
\qquad x\in\R^d.
\]
\end{proposition}

\begin{proof}
If such a function $h$ existed, applying the heat semigroup and using Dynkin's formula would give
\begin{equation}\label{pp1}
(P_t h)(x)-h(x)
+b_*\int_0^t(P_s(qh))(x)\,ds
\ge a\int_0^t(P_s q)(x)\,ds.
\end{equation}
Since $h\in C_0$, we have $\|P_t h\|_\infty\to0$. Also,
$G|qh|\le\|h\|_\infty Gq<\infty$. Letting $t\to\infty$ in \eqref{pp1}, we
obtain
\begin{equation}\label{eq:auxiliary-endpoint-potential}
 h(x) \le b_*(Th)(x) -a(Gq)(x).
\end{equation}

By Lemma~\ref{lem:spectral} and the Krein--Rutman theorem
(see \cite{Schaefer}), there exists a nonzero positive functional
$\ell\in C_0(\R^d)^*$ satisfying $T^*\ell=b_*^{-1}\ell$.
By the definition of the dual operator,
\[
\ell(Tf)=(T^*\ell)(f)=b_*^{-1}\ell(f),
\qquad f\in C_0(\R^d).
\]
The Riesz representation theorem gives a nonzero finite positive
Radon measure $\mu_\ell$ such that
\[
\ell(f)=\int_{\R^d}f(x)\,\mu_\ell(dx),
\qquad f\in C_0(\R^d).
\]
Since $(Gq)(x)>0$ for every $x\in\R^d$, we have
\[
\ell(Gq)=\int_{\R^d}(Gq)(x)\,\mu_\ell(dx)>0.
\]
Applying the positive functional $\ell$ to
\eqref{eq:auxiliary-endpoint-potential}, we obtain
\[
\begin{aligned}
	\int_{\R^d}h(x)\,\mu_\ell(dx)
	&\le b_*\int_{\R^d}(Th)(x)\,\mu_\ell(dx)
	-a\int_{\R^d}(Gq)(x)\,\mu_\ell(dx)\\
	&=\int_{\R^d}h(x)\,\mu_\ell(dx)
	-a\int_{\R^d}(Gq)(x)\,\mu_\ell(dx)\\
	&<\int_{\R^d}h(x)\,\mu_\ell(dx).
\end{aligned}
\]
This contradiction proves the proposition.
\end{proof}

To prove uniform integrability at $b_*$ and for nearby $\kappa>b_*$,
we compare the size-biased partition function with a functional of a
single Brownian motion.

\begin{proposition}\label{prop:scalar}
	Fix $\kappa\ge0$. Let $\eta$ be a standard real Brownian motion,
	independent of $D$. For $t\ge0$, define
	\[
	F_t^\kappa(\eta)
	=E_0^D\!\left[
	\exp\left\{
	\sqrt\kappa\int_0^t q(D_s)\,d\eta_s
	+\kappa\int_0^t q(D_s)\,ds
	-\frac{\kappa}{2}\int_0^t q(D_s)^2\,ds
	\right\}
	\right].
	\]
	If $\{F_n^\kappa(\eta):n\in\mathbb N\}$ is tight under the law
	of $\eta$, then $(Z_t^\kappa(0))_{t\ge0}$ is uniformly integrable. Consequently, weak disorder holds at $\kappa$.
\end{proposition}

\begin{proof}
	
	Fix $t\ge0$ and a bounded measurable function
	$\Phi:\R_+\to\R$. Write $Z_t^\kappa(0;W)$ to display its
	dependence on the environment. By \eqref{eq:polymer-definition}
	and Fubini's theorem,
	\[
	\E\!\left[
	Z_t^\kappa(0;W)\Phi(Z_t^\kappa(0;W))
	\right]
	=
	E_0^B\E\!\left[
	e^{\sqrt\kappa H_t(B)-\kappa t/2}
	\Phi(Z_t^\kappa(0;W))
	\right].
	\]
	Define the shifted environment on $[0,t]$ by
	\begin{equation}\label{lq}
	W^{\kappa,B}(ds,x)
	=W(ds,x)+\sqrt\kappa\,q(x-B_s)\,ds.
	\end{equation}
	For each fixed path $B$, the Cameron--Martin formula gives
	\[
	\E\!\left[
	e^{\sqrt\kappa H_t(B)-\kappa t/2}
	\Phi(Z_t^\kappa(0;W))
	\right]
	=
	\E\!\left[
	\Phi(Z_t^\kappa(0;W^{\kappa,B}))
	\right].
	\]
	Write $\widehat\Pp=P_0^B\otimes\Pp^W$ for the product law of $B$ and
	the environment, and $\widehat\E$ for its expectation. Adjoin an
	independent Brownian motion $\bar B$ started at $0$, retaining this
	notation for the extended law. Substituting \eqref{lq}
	into \eqref{eq:polymer-definition}, we obtain
	\begin{equation}\label{oo1}
	Z_t^\kappa(0;W^{\kappa,B})
	=
	E_0^{\bar B}\!\left[
	\exp\left\{
	\sqrt\kappa\int_0^t W(ds,\bar B_s)-\frac{\kappa t}{2}
	+\kappa\int_0^t q(\bar B_s-B_s)\,ds
	\right\}
	\right]
	=:\widetilde Z_t^\kappa.
	\end{equation}
	Combining the preceding identities yields
	\begin{equation}\label{m2}
	\E\!\left[
	Z_t^\kappa(0)\Phi(Z_t^\kappa(0))
	\right]
	=
	\widehat\E\!\left[\Phi(\widetilde Z_t^\kappa)\right].
	\end{equation}

Define 
$d\eta_s=W(ds,B_s)$. Conditionally on $(B,\bar B)$, $\eta$ has the law of a standard Brownian motion, since $q(0)=1$ and both paths are independent of $W$. Thus $\eta$ is independent of the pair $(B,\bar B)$. For fixed paths, Gaussian orthogonal projection gives
\[
 W(ds,x)=q(x-B_s)d\eta_s+dW_B^\perp(s,x),
\]
and we set
\[
R_t(B,\bar B):=\int_0^t dW_B^\perp(s,\bar B_s).
\]
Given $(B,\bar B)$, the residual $R_t(B,\bar B)$ is jointly Gaussian with $\eta$ and is orthogonal to $\int_0^t g(s)\,d\eta_s$ for every deterministic $g\in L^2[0,t]$. Hence, conditionally on $(\eta,B,\bar B)$, it is centered Gaussian with variance
\[
\int_0^t\bigl[1-q(\bar B_s-B_s)^2\bigr]\,ds.
\]
Therefore, 
\begin{equation}\label{w1}
\begin{aligned}
 \widehat\E[e^{\sqrt\kappa R_t-\kappa t/2}
                       \mid\eta,B,\bar B]
 &=\exp\left\{\frac\kappa2\int_0^t
                 [1-q(\bar B_s-B_s)^2]ds-\frac{\kappa t}2\right\}\\
 &=\exp\left\{-\frac\kappa2\int_0^t
                            q(\bar B_s-B_s)^2ds\right\}.
\end{aligned}
\end{equation}
Combining \eqref{oo1} and \eqref{w1} gives
\[
\begin{aligned}
 \widehat\E[\widetilde Z_t^\kappa\mid\eta,B]
 &=E_0^{\bar B}\exp\biggl\{\sqrt\kappa\int_0^t q(\bar B_s-B_s)\,d\eta_s\\
 &\hspace{3em}+\kappa\int_0^t q(\bar B_s-B_s)\,ds
 -\frac{\kappa}{2}\int_0^t q(\bar B_s-B_s)^2\,ds\biggr\}.
\end{aligned}
\]
Averaging over $B$ conditional on $\eta$, we see that the difference
${\bar B}-B$ has the law of a Brownian motion with generator $\Delta$, independent of $\eta$. Consequently,
$\widehat\E[\widetilde Z_t^\kappa\mid\eta]=F_t^\kappa(\eta)$.

 For any $L,K>0$, conditional Markov's inequality gives
\begin{equation}\label{p1}
\begin{aligned}
 \widehat\Pp(\widetilde Z_t^\kappa>K)
 &\le\Pp(F_t^\kappa>L)
     +\widehat\E[\mathbf1_{\{F_t^\kappa\le L\}}
                    \widehat\Pp(\widetilde Z_t^\kappa>K\mid\eta)]\\
 &\le\Pp(F_t^\kappa>L)
       +K^{-1}\E[\mathbf1_{\{F_t^\kappa\le L\}}F_t^\kappa]
 \le\Pp(F_t^\kappa>L)+L/K.
\end{aligned}
\end{equation}
Combining \eqref{m2} and \eqref{p1} with the tightness of
$(F_n^\kappa)_{n\ge1}$ gives
\[
\lim_{K\to\infty}\sup_{n\ge1}
\E\!\left[
Z_n^\kappa(0)\mathbf1_{\{Z_n^\kappa(0)>K\}}
\right]
=
\lim_{K\to\infty}\sup_{n\ge1}
\widehat\Pp(\widetilde Z_n^\kappa>K)
=0.
\]
Thus $(Z_n^\kappa(0))_{n\ge1}$ is uniformly integrable and
converges to $Z_\infty^\kappa(0)$ in $L^1$.

For each $t\ge0$, passing to the $L^1$ limit in
\[
Z_t^\kappa(0)
=\E[Z_n^\kappa(0)\mid\mathcal F_t^W],
\qquad n\in\mathbb N,\quad n\ge t,
\]
yields
\[
Z_t^\kappa(0)
=\E[Z_\infty^\kappa(0)\mid\mathcal F_t^W].
\]
 Hence
$(Z_t^\kappa(0))_{t\ge0}$ is uniformly integrable, and weak
disorder holds at $\kappa$.
\end{proof}

\subsection{A bounded gauge and a local limit after rank-one deletion}

We next subtract a small rank-one kernel from the critical
Feynman--Kac kernel, construct a bounded positive harmonic
function for the resulting kernel, and establish a local
limit theorem for its iterates.

Set $b=b_*$, $P=P_1$, and
$S_t=e^{t(\Delta+bq)}$, with $S=S_1$.
For the bounded function $\psi$ in \eqref{oo}, the Feynman--Kac formula gives
\begin{equation}\label{ll}
	(S_t\psi)(x)
	=
	E_x^D\!\left[
	e^{b\int_0^tq(D_s)\,ds}\psi(D_t)
	\right]
	=\psi(x),
	\qquad t\ge0,\quad x\in\R^d.
\end{equation}
Let $0\ne f\in C_c^\infty(\R^d)$ be nonnegative.
Fix an open ball $C$ containing $\operatorname{supp}f$. For $\lambda\ge0$, take the following version of the time-one Feynman--Kac density restricted to paths remaining in $C$:
\begin{equation}\label{ebridge}
 S_\lambda^C(x,y):=\rho_1(y-x)E_{x,y}^{\mathrm{bridge}}\!\left[
 e^{\lambda\int_0^1q(D_s)\,ds}
 \mathbf1_{\{D_s\in C\ \text{for all }s\in[0,1]\}}\right],
 \qquad x,y\in\R^d .
\end{equation}
Here $E_{x,y}^{\mathrm{bridge}}$ denotes expectation for the unit-time bridge of the difference motion $D$ from $x$ to $y$. Conditioning on the endpoint shows that $S_\lambda^C(x,y)\,dy$ is the weighted endpoint measure for paths started at $x$ and remaining in $C$. Write $S_C=S_b^C$.

Define the operators
\[
A:=S-P-\varepsilon f\otimes f,
\]
and
\begin{equation}\label{sy1}
	\mathcal B:=S-\varepsilon f\otimes f=P+A.
\end{equation}
The integral kernels corresponding to the operators $A$ and
$\mathcal B$ are denoted by $A(x,y)$ and $\mathcal B(x,y)$,
respectively. In particular,
\[
(\mathcal Bv)(x)
=\int_{\R^d}\mathcal B(x,y)v(y)\,dy,
\qquad
\mathcal B(x,y)=S(x,y)-\varepsilon f(x)f(y).
\]

\begin{lemma}\label{lem:rank-one-gauge}

	There exists $\varepsilon_0=\varepsilon_0(f,C)>0$ such that,
	for every $0<\varepsilon<\varepsilon_0$, 
    the following statements hold.
	\begin{enumerate}
		\item
		\[
		\rho_1(x-y)\le\mathcal B(x,y)\le e^b\rho_1(x-y),
		\qquad x,y\in\R^d.
		\]
		
		\item There exists $H\in1+C_0(\R^d)$ such that
		\[
		H(x)\ge1,\qquad
		(\mathcal BH)(x)=H(x),\qquad
		(\mathcal B^n1)(x)\uparrow H(x),
		\qquad x\in\R^d.
		\]
		
		\item
		\[
		\langle g,\mathcal Bg\rangle_{L^2}
		\ge-\frac12\|g\|_{L^2}^2,
		\qquad g\in L^2(\R^d).
		\]
		
		\item
		\[
		\varepsilon f(x)f(y)\le S_C(x,y),
		\qquad x,y\in\R^d.
		\]
	\end{enumerate}
\end{lemma}

\begin{proof}
	
	We first prove \textup{(1)}. Let $R:=S-P$. The Brownian-bridge
	representation gives
	\[
	R(x,y)
	=
	\rho_1(y-x)E_{x,y}^{\mathrm{bridge}}\!\left[
	e^{b\int_0^1q(D_s)\,ds}-1
	\right].
	\]
    It follows that $R(x,y)$ is continuous and strictly positive.
	Let $K=\operatorname{supp}f$. By compactness,
	$
	\min_{(x,y)\in K\times K}R(x,y)>0.
	$ 
	Choose $\varepsilon>0$ sufficiently small that
	$$
	\varepsilon\|f\|_\infty^2
	\le \frac12\min_{(x,y)\in K\times K}R(x,y).
	$$
	Since $f(x)f(y)=0$ outside $K\times K$, this implies
	$
	0\le\varepsilon f(x)f(y)\le\frac12R(x,y),
x,y\in\R^d.
	$
	Consequently,
	\[
	\begin{aligned}
		\mathcal B(x,y)
		&=\rho_1(y-x)+R(x,y)-\varepsilon f(x)f(y)\\
		&\ge \rho_1(y-x)+\frac12R(x,y)
		\ge \rho_1(y-x).
	\end{aligned}
	\]
	On the other hand, since $\varepsilon f(x)f(y)\ge0$ and
	$0\le q\le1$,
	\[
	\begin{aligned}
		\mathcal B(x,y)
		&\le S(x,y)
		=\rho_1(y-x)E_{x,y}^{\mathrm{bridge}}\!\left[
		e^{b\int_0^1q(D_s)\,ds}
		\right]
		\le e^b\rho_1(y-x).
	\end{aligned}
	\]
	These two bounds prove \textup{(1)}.

We now prove \textup{(2)}. Whenever the series converges
absolutely, write
$
(\mathcal Gv)(x):=\sum_{j=0}^\infty(P^jv)(x),
$
and set $\mathcal T=\mathcal GR$. For bounded measurable $v$, Duhamel's formula gives
\[
(Rv)(x)
=
b\int_0^1
\bigl[P_s(qS_{1-s}v)\bigr](x)\,ds.
\]
Since $\|S_sv\|_\infty\le e^{bs}\|v\|_\infty$, we obtain
\begin{equation}\label{eq:discrete-greenbound}
	\begin{aligned}
		|(\mathcal Tv)(x)|
		&\le\sum_{j=0}^\infty(P^j|Rv|)(x)
		\le be^b\|v\|_\infty
		\sum_{j=0}^\infty\int_0^1(P_{j+s}q)(x)\,ds
		=be^b\|v\|_\infty(Gq)(x).
	\end{aligned}
\end{equation}
The proof of Lemma~\ref{lem:spectral} gives $Gq\in C_0(\R^d)$.
Thus the series defining $\mathcal Tv$ converges absolutely
and $\|\mathcal Tv\|_\infty\le C\|v\|_\infty$.

The heat-kernel estimate $\|\nabla\rho_s\|_1\le Cs^{-1/2}$
implies
\[
\|\nabla Rv\|_\infty
\le be^b\|v\|_\infty
\int_0^1\|\nabla\rho_s\|_1\,ds
\le C\|v\|_\infty.
\]
Absolute convergence permits applying $P$ term by term, so
$\mathcal Tv=Rv+P\mathcal Tv$. Consequently,
\[
\|\nabla\mathcal Tv\|_\infty
\le
\|\nabla Rv\|_\infty
+\|\nabla\rho_1\|_1\|\mathcal Tv\|_\infty
\le C\|v\|_\infty.
\]
Together with \eqref{eq:discrete-greenbound}, this shows that
$\mathcal T$ maps bounded measurable functions into $C_0(\R^d)$.
The image of the unit ball is equicontinuous and vanishes
uniformly at infinity. The Arzel\`a--Ascoli theorem therefore implies
that $\mathcal T$ is compact on $C_0(\R^d)$.

For the bounded function $\psi$ in \eqref{oo}, \eqref{ll} gives
$R\psi=\psi-P\psi$. Since $\psi\in C_0(\R^d)$,
$\|P^N\psi\|_\infty\to0$, and
\[
\sum_{j=0}^{N-1}(P^jR\psi)(x)
=
\psi(x)-(P^N\psi)(x)
\longrightarrow\psi(x).
\]
Thus $\mathcal T\psi=\psi$.

Since $f$ is compactly supported, the heat-kernel bound gives
\[
\|P^jf\|_\infty
\le (4\pi j)^{-d/2}\|f\|_1,
\qquad j\ge1.
\]
Since $d\ge3$,
\[
\sum_{j=0}^\infty\|P^jf\|_\infty
\le \|f\|_\infty
+(4\pi)^{-d/2}\|f\|_1
\sum_{j=1}^\infty j^{-d/2}<\infty.
\]
Thus the series $\mathcal Gf=\sum_{j=0}^\infty P^jf$
converges absolutely and uniformly.
Since each $P^jf$ belongs to $C_0(\R^d)$, so does $\mathcal Gf$.
Moreover, $f\ge0$, $f\not\equiv0$, and the strict positivity
of the heat kernel imply
\[
(\mathcal Gf)(x)\ge(Pf)(x)>0,
\qquad x\in\R^d.
\]
Define
\[
\mathcal T_\varepsilon
:=
\mathcal GA
=
\mathcal T-\varepsilon(\mathcal Gf)\otimes f.
\]
This operator is compact on $C_0(\R^d)$ and maps bounded
measurable functions into $C_0(\R^d)$.
The estimate obtained in the proof of \textup{(1)} gives
\[
\frac12R\le A\le R,
\qquad
\frac12\mathcal T\le\mathcal T_\varepsilon\le\mathcal T.
\]
In particular, $\mathcal T_\varepsilon$ is positive.
Since $\mathcal T\psi=\psi$, positivity and the
spectral-radius formula yield
\[
\mathcal T_\varepsilon^n\psi\ge2^{-n}\psi,
\qquad
r(\mathcal T_\varepsilon)\ge\frac12.
\]

Applying the argument from the proof of
Proposition~\ref{prop:auxiliary-endpoint-obstruction} to
$\mathcal T_\varepsilon$, we obtain a nonzero
positive functional $\ell\in C_0(\R^d)^*$ satisfying
$
\mathcal T_\varepsilon^*\ell
=r(\mathcal T_\varepsilon)\ell.
$
The same positivity argument gives
$\ell(\psi)>0$ and $\ell(\mathcal Gf)>0$.
We obtain
\[
\begin{aligned}
	r(\mathcal T_\varepsilon)\ell(\psi)
	&=\ell(\mathcal T_\varepsilon\psi)
	=\ell(\psi)
	-\varepsilon\ell(\mathcal Gf)
	\int_{\R^d}f(x)\psi(x)\,dx
	<\ell(\psi).
\end{aligned}
\]
Dividing by $\ell(\psi)>0$ yields
\begin{equation}\label{eq:strict-deletionradius}
	r(\mathcal T_\varepsilon)<1.
\end{equation}

Since $\mathcal T_\varepsilon1\in C_0(\R^d)$,
\eqref{eq:strict-deletionradius} ensures uniform convergence
of the nonnegative series
$
H(x)
:=
1+\sum_{j=0}^\infty
\bigl[\mathcal T_\varepsilon^j
(\mathcal T_\varepsilon1)\bigr](x).
$
Thus $H\in1+C_0(\R^d)$, $H\ge1$, and
\[
H(x)-1
=(\mathcal T_\varepsilon H)(x)
=\sum_{j=0}^\infty(P^jAH)(x).
\]
Applying $P$ term by term gives
\[
\begin{aligned}
	\bigl[P(H-1)\bigr](x)
	&=H(x)-1-(AH)(x),
\end{aligned}
\]
and hence $(\mathcal BH)(x)
=1+\bigl[P(H-1)\bigr](x)+(AH)(x)
=H(x).$
Finally, $\mathcal B1\ge P1=1$ and $\mathcal BH=H$ imply
\[
1\le(\mathcal B^n1)(x)
\le(\mathcal B^{n+1}1)(x)\le H(x).
\]
Let $h(x):=\lim_{n\to\infty}(\mathcal B^n1)(x)$.
Monotone convergence gives $\mathcal Bh=h$. Moreover,
\[
0\le P^N(h-1)\le P^N(H-1),
\qquad
\|P^N(H-1)\|_\infty\longrightarrow0.
\]
Iterating $h=Ph+Ah$ therefore yields
\[
h(x)
=(P^Nh)(x)+\sum_{j=0}^{N-1}(P^jAh)(x)
\longrightarrow1+(\mathcal T_\varepsilon h)(x).
\]
Hence $h-1=\mathcal T_\varepsilon h\in C_0(\R^d)$.
Subtracting $H-1=\mathcal T_\varepsilon H$ gives
\[
(I-\mathcal T_\varepsilon)(h-H)=0.
\]
By \eqref{eq:strict-deletionradius},
$I-\mathcal T_\varepsilon$ is invertible on $C_0(\R^d)$,
so $h=H$. Thus $(\mathcal B^n1)(x)\uparrow H(x)$ for every
$x\in\R^d$, proving \textup{(2)}.

We now prove \textup{(3)}. Since $S=e^{\Delta+bq}$ is self-adjoint and positive semidefinite on
$L^2(dx)$, we have
\[
 \langle v,\mathcal Bv\rangle
 =\langle v,Sv\rangle-\varepsilon|\langle f,v\rangle|^2
 \ge-\varepsilon\|f\|_2^2\|v\|_2^2,
 \qquad v\in L^2(dx).
\]
Reducing $\varepsilon_0$ so that
$\varepsilon_0\|f\|_2^2\le1/2$, we obtain
\[
\langle v,\mathcal Bv\rangle
\ge-\frac12\|v\|_2^2,
\qquad v\in L^2(dx),
\]
which proves \textup{(3)}.

For \textup{(4)}, choose a concentric ball $C'$ with $K\subset C'\Subset C$ and set $r=\operatorname{dist}(C',\partial C)>0$. For $x,y\in K$, the line segment from $x$ to $y$ lies in $C'$. The bridge of $D$ from $x$ to $y$ is this linear path plus a centered bridge of duration one. The latter stays within distance $r$ of zero with probability $p_r>0$, independently of the endpoints. Since the Feynman--Kac weight is at least one,
\[
 \inf_{x,y\in K}S_C(x,y)\ge p_r\inf_{x,y\in K}\rho_1(y-x)>0.
\]
Reduce $\varepsilon_0$ so that $\varepsilon_0\|f\|_\infty^2$ is at most this lower bound. Outside $K\times K$ the product $f(x)f(y)$ vanishes, which proves \textup{(4)}.
\end{proof}

Take $\varepsilon>0$ sufficiently small
for all conclusions of Lemma~\ref{lem:rank-one-gauge} to hold.
Lemma~\ref{lem:regen-local-limit} then applies to the
resulting kernel $\mathcal B$ and the function $H$.

\begin{lemma}\label{lem:regen-local-limit}
	Let $\mathcal B$ and $H$ be as in
    Lemma~\ref{lem:rank-one-gauge}. Then, for fixed
    $x,y\in\R^d$,
	\begin{equation}\label{eq:regen-local-limit}
		\mathcal B^n(x,y)\sim(4\pi n)^{-d/2}H(x)H(y),
		\qquad n\to\infty.
	\end{equation}
	Moreover, there exists $C_1>0$ such that
	\[
	\mathcal B^n(x,y)\le C_1n^{-d/2},
	\qquad n\ge1,\quad x,y\in\R^d.
	\]
\end{lemma}

\begin{proof}
	Define
	\[
	\mathsf P(x,dy)
	:=\mathcal B(x,y)\frac{H(y)}{H(x)}\,dy, 
	\qquad
	m_H(dx):=H(x)^2\,dx,
	\]
	and
	\[
	(\mathsf Pv)(x)
	:=\int_{\R^d}v(y)\,\mathsf P(x,dy).
	\]
	By Lemma~\ref{lem:rank-one-gauge} \textup{(2)},
	$\mathcal BH=H$, so $\mathsf P1=1$. The operator $\mathsf P$ is a self-adjoint contraction on
    $L^2(m_H)$. Indeed, symmetry of $\mathcal B$ gives
	\[
	m_H(dx)\,\mathsf P(x,dy)
	=
	m_H(dy)\,\mathsf P(y,dx).
	\]
	This reversibility identity gives symmetry on $L^2(m_H)$.
	Together with $\mathsf P1=1$, it also implies that $m_H$
	is invariant. Jensen's inequality then yields the
	$L^2(m_H)$ contraction property. 
	
For $v\in L^2(m_H)$, we have $Hv\in L^2(dx)$.
The definition of $\mathsf P$ and
Lemma~\ref{lem:rank-one-gauge}\textup{(3)} give
\[
\begin{aligned}
	\langle v,\mathsf Pv\rangle_{L^2(m_H)}
	&=\langle Hv,\mathcal B(Hv)\rangle_{L^2(dx)}
	\ge-\frac12\|Hv\|_{L^2(dx)}^2
	=-\frac12\|v\|_{L^2(m_H)}^2.
\end{aligned}
\]
This yields
$
\sigma(\mathsf P)\subseteq[-1/2,1].
$ For each $x\in\R^d$, the transition measure $\mathsf P(x,dy)$
has density $\mathcal B(x,y)/(H(x)H(y))$ with respect to
$m_H(dy)$. This density is uniformly bounded in $x,y$
by Lemma~\ref{lem:rank-one-gauge}\textup{(1)} and $H\ge1$.
	
	For $v\in L^1(m_H)\cap L^2(m_H)$, set
	$
	\mathcal E(v)
	:=\langle v,(I-\mathsf P)v\rangle_{m_H}.
	$
    Then
	\begin{align*}
		\mathcal E(v)
		&=\frac12\iint
		\mathcal B(x,y)H(x)H(y)|v(x)-v(y)|^2\,dx\,dy\\
		&\ge\frac12\iint
		\rho_1(x-y)|v(x)-v(y)|^2\,dx\,dy\\
		&=\langle v,(I-P)v\rangle_{L^2(dx)}=\int_{\R^d}
		(1-e^{-|\xi|^2})|\widehat v(\xi)|^2\,d\xi,
	\end{align*}
    where $\widehat v$ denotes the unitary Fourier transform of $v$.
    For $0<r\le1$, we have
	\begin{equation}\label{de}
	\|v\|_{L^2(dx)}^2
	\le Cr^d\|v\|_{L^1(dx)}^2
	+Cr^{-2}\mathcal E(v).
	\end{equation}
    Set $A_v:=\|v\|_{L^2(m_H)}^2$ and $B_v:=\|v\|_{L^1(m_H)}^2$.
    Since $1\le H\le\|H\|_\infty$, \eqref{de} gives
	\[
	A_v\le Cr^dB_v+Cr^{-2}\mathcal E(v).
	\]
    Fix $M>0$. If $v\ne0$ and $A_v\le MB_v$, choose
	$
	r^d=\min\left\{1,\frac{A_v}{2CB_v}\right\}.
	$ 
	Then $Cr^dB_v\le A_v/2$, and therefore
	\begin{equation}\label{eq:restricted-nash-calculation}
		\mathcal E(v)
		\ge cr^2A_v
		\ge c_M A_v^{1+2/d}B_v^{-2/d},
	\end{equation}
	where
	$
	c:=\frac{1}{2C},
	c_M:=c\min\left\{(2C)^{-2/d},\,M^{-2/d}\right\}.
	$
	
	Let $u\ge0$ belong to $L^1(m_H)$ with
	$\|u\|_{L^1(m_H)}=1$, and set $v_n:=\mathsf P^nu$.
	Since $m_H$ is invariant and $v_n\ge0$, we have
	\[
	\|v_n\|_{L^1(m_H)}
	=\int_{\R^d}\mathsf P^nu\,dm_H
	=\int_{\R^d}u\,dm_H=1.
	\]
    Since
    $\|v_1\|_\infty\le C\|u\|_{L^1(m_H)}=C$, we have
	\[
	\|v_n\|_\infty
	=\|\mathsf P^{n-1}v_1\|_\infty
	\le\|v_1\|_\infty\le C,
	\qquad n\ge1.
	\]
	Consequently,
	\[
	a_n:=\|v_n\|_{L^2(m_H)}^2
	\le\|v_n\|_\infty\|v_n\|_{L^1(m_H)}
	\le C.
	\]
    Since $A_{v_n}=a_n\le C$ and $B_{v_n}=1$,
	\eqref{eq:restricted-nash-calculation} applies uniformly
	in $n\ge1$ with a fixed $M\ge C$.
	Since
	\[
	1-\lambda^2\ge\frac12(1-\lambda),
	\qquad \lambda\in[-1/2,1],
	\]
	the spectral theorem yields
	$I-\mathsf P^2\ge\frac12(I-\mathsf P)$. Consequently,
	\[
	\begin{aligned}
		a_n-a_{n+1}
		&=\langle v_n,(I-\mathsf P^2)v_n\rangle_{m_H}
		\ge\frac12\mathcal E(v_n)
		\ge \tfrac12 c_M a_n^{1+2/d},
		\qquad n\ge1.
	\end{aligned}
	\]
	Since $\|v_n\|_{L^1(m_H)}=1$, we have $a_n>0$.
	The mean value theorem now gives
	$$
	a_{n+1}^{-2/d}-a_n^{-2/d}
	\ge\frac{c_M}{d}.
	$$
    Together with $a_1\le C$, this implies $a_n\le Cn^{-d/2}$. For $0\ne u\in L^1(m_H)$, positivity gives
	$
	|\mathsf P^n u|\le \mathsf P^n|u|.
	$
	Applying the preceding estimate to
	$|u|/\|u\|_{L^1(m_H)}$ yields
	\begin{equation}\label{tn}
	\|\mathsf P^n u\|_{L^2(m_H)}^2
	\le Cn^{-d/2}\|u\|_{L^1(m_H)}^2.
	\end{equation}
Write $\|\cdot\|_{p\to q}$ for the operator norm from
$L^p(m_H)$ to $L^q(m_H)$. Taking the supremum in \eqref{tn} over $\|u\|_{L^1(m_H)}=1$, we obtain
	\[
	\|\mathsf P^n\|_{1\to 2}^2
	\le Cn^{-d/2},
	\qquad n\ge1.
	\]
For $m\ge1$, self-adjointness and composition give
	\[
	\begin{aligned}
		\|\mathsf P^{2m}\|_{1\to\infty}
		&\le
		\|\mathsf P^m\|_{2\to\infty}
		\|\mathsf P^m\|_{1\to2}
		=\|\mathsf P^m\|_{1\to2}^2
		\le Cm^{-d/2}.
	\end{aligned}
	\]
For $m\ge1$, the $L^\infty(m_H)$ contraction property gives
\[
\|\mathsf P^{2m+1}\|_{1\to\infty}
\le \|\mathsf P^{2m}\|_{1\to\infty}
\le Cm^{-d/2}.
\]
Since $2m+1\le3m$, this yields the required bound for odd
$n\ge3$. The bounded one-step density gives
$\|\mathsf P\|_{1\to\infty}\le C$.
Therefore
$
\|\mathsf P^n\|_{1\to\infty}
\le Cn^{-d/2},
 n\ge1.
$
	
	For $m_H\otimes m_H$-almost every $(x,y)$, the density of $\mathsf P^n$ satisfies
	\[
	\frac{\mathcal B^n(x,y)}{H(x)H(y)}
	\le \|\mathsf P^n\|_{1\to\infty}
	\le Cn^{-d/2}.
	\]
	For each fixed $n$, the kernel $\mathcal B^n(x,y)$ is jointly continuous by iterated convolution and dominated convergence, using the continuity of $\mathcal B$ and the Gaussian bound $\mathcal B(x,y)\le e^b\rho_1(x-y)$ from Lemma~\ref{lem:rank-one-gauge}. Since $H$ is continuous and $H\ge1$, the density on the left is continuous and $m_H\otimes m_H$ has full support. The bound therefore holds for every $(x,y)$. As $H$ is bounded, it follows that
	\begin{equation}\label{uu1}
	\mathcal B^n(x,y)\le Cn^{-d/2},
	\qquad n\ge1,\quad x,y\in\R^d.
	\end{equation}

    We next prove \eqref{eq:regen-local-limit}. Recall from
    Lemma~\ref{lem:rank-one-gauge} that $A=\mathcal B-P$.
	Since $A1=\mathcal B1-1$, we get
	\[
	\begin{aligned}
		\sum_{j=0}^{N-1}(\mathcal B^jA1)(x)
		&=\sum_{j=0}^{N-1}
		\bigl[(\mathcal B^{j+1}1)(x)-(\mathcal B^j1)(x)\bigr]
		=(\mathcal B^N1)(x)-1.
	\end{aligned}
	\]
	By Lemma~\ref{lem:rank-one-gauge}\textup{(2)},
	letting $N\to\infty$ yields
	\begin{equation}\label{eq:regen-telescope}
		\sum_{j=0}^{\infty}(\mathcal B^jA1)(x)=H(x)-1.
	\end{equation}
	In particular,
	\[
	\sum_{j=m}^{\infty}(\mathcal B^jA1)(x)
	=H(x)-(\mathcal B^m1)(x)\longrightarrow0.
	\]
	For integers $m\ge1$ and $n>2m$, define
	$
	E_{n,m}:=\mathcal B^mP^{n-2m}\mathcal B^m.
	$
	Since $A\ge0$, positivity gives
	$P^j\le\mathcal B^j$ and hence $E_{n,m}\le\mathcal B^n$.
	The telescoping identity
	\[
	\mathcal B^N-P^N
	=\sum_{i=0}^{N-1}\mathcal B^iAP^{N-i-1}
	\]
	with $N=n-2m$ gives
	\[
	\begin{aligned}
		\mathcal B^n-E_{n,m}
		&=\sum_{i=0}^{n-2m-1}
		\mathcal B^{m+i}AP^{n-2m-i-1}\mathcal B^m
		\le
		\sum_{j=m}^{n-m-1}
		\mathcal B^jA\mathcal B^{n-j-1}.
	\end{aligned}
	\]
    For $j\le(n-1)/2$, \eqref{uu1} gives
	\[
	\begin{aligned}
		(\mathcal B^jA\mathcal B^{n-j-1})(x,y)
		&=\iint
		\mathcal B^j(x,z)A(z,w)
		\mathcal B^{n-j-1}(w,y)\,dz\,dw
		\le Cn^{-d/2}(\mathcal B^jA1)(x).
	\end{aligned}
	\]
	For $j>(n-1)/2$, symmetry gives the bound
	$Cn^{-d/2}(\mathcal B^{n-j-1}A1)(y)$.
	Thus \eqref{eq:regen-telescope} implies
	\begin{equation}\label{qa}
	\begin{aligned}
		0&\le n^{d/2}
		\bigl(\mathcal B^n-E_{n,m}\bigr)(x,y)
		\le C\bigl[
		H(x)-(\mathcal B^m1)(x)
		+H(y)-(\mathcal B^m1)(y)
		\bigr].
	\end{aligned}
	\end{equation}
	The right-hand side is independent of $n$ and tends to zero
	as $m\to\infty$ by
	Lemma~\ref{lem:rank-one-gauge}\textup{(2)}.
	
	For fixed $m$, the heat-kernel formula gives
	\[
	\begin{aligned}
		\frac{E_{n,m}(x,y)}{(4\pi n)^{-d/2}}
		&=\iint
		\mathcal B^m(x,z)\mathcal B^m(w,y)
		\left(\frac{n}{n-2m}\right)^{d/2}
		e^{-|z-w|^2/[4(n-2m)]}\,dz\,dw.
	\end{aligned}
	\]
    By symmetry,
	\[
	\begin{aligned}
		\iint\mathcal B^m(x,z)\mathcal B^m(w,y)\,dz\,dw
		&=(\mathcal B^m1)(x)(\mathcal B^m1)(y)
		\le H(x)H(y)
		\le \|H\|_\infty^2<\infty.
	\end{aligned}
	\]
    Dominated convergence therefore yields
	\[
	\lim_{n\to\infty}
	\frac{E_{n,m}(x,y)}{(4\pi n)^{-d/2}}
	=(\mathcal B^m1)(x)(\mathcal B^m1)(y).
	\]
	Combining this limit with \eqref{qa} and then
    letting $m\to\infty$, we obtain
	\[
	\lim_{n\to\infty}(4\pi n)^{d/2}\mathcal B^n(x,y)
	=H(x)H(y).
	\]
This proves \eqref{eq:regen-local-limit}.
\end{proof}

For $\delta\ge0$, keep the rank-one kernel $\varepsilon f\otimes f$
fixed and set
\[
 \mathcal B_\delta=e^{\Delta+(b+\delta)q}-\varepsilon f\otimes f,
 \qquad \mathcal B_0=\mathcal B .
\]
The following estimates control the effect of increasing the
environmental amplitude on the iterates of this kernel.

\begin{lemma}\label{lem:gap-perturb}
There exist $\delta_0>0$ and $C<\infty$ such that, for
$0\le\delta\le\delta_0$, $n\ge1$, and $x,y\in\R^d$,
\begin{align}
 (\mathcal B_\delta^n1)(x)&\le 2H(x)-1,\notag\\
 0\le\mathcal B_\delta^n(x,y)-\mathcal B^n(x,y)
   &\le C\delta n^{-d/2}. \label{eq:gap-kernel-stability}
\end{align}
\end{lemma}

\begin{proof}
Let $E_\delta=\mathcal B_\delta-\mathcal B$. For
$u=\int_0^1 q(D_s)\,ds\in[0,1]$ and $0\le\delta\le1$,
\[
 0\le e^{bu}(e^{\delta u}-1)
 \le e^{b+1}\delta u
 \le (e^{b+1}/b)\delta(e^{bu}-1).
\]
Applying the Brownian bridge representation of the time-one kernels
gives $0\le E_\delta(x,y)\le(e^{b+1}/b)\delta(S-P)(x,y)$.
By the choice of $\varepsilon$ in the proof of
Lemma~\ref{lem:rank-one-gauge}, $S-P\le2A$ pointwise.
Thus $0\le E_\delta\le c\delta A$ as kernels and
$\mathcal B_\delta=\mathcal B+E_\delta\ge\mathcal B\ge P$.

Put $w=2H-1\ge1$. Since $\mathcal BH=H$ and $P1=1$,
\[
 (I-\mathcal B)w=A1,\qquad
 E_\delta w\le c\delta\|w\|_\infty A1 .
\]
Choose $\delta_0\le1$ such that $c\delta_0\|w\|_\infty\le1/2$.
Then $(I-\mathcal B_\delta)w\ge A1/2$, so positivity gives
$\mathcal B_\delta^m1\le\mathcal B_\delta^m w\le w$ for $m\ge0$.
Moreover, telescoping yields, for $m\ge1$,
\[
 \begin{aligned}
 \sum_{j=0}^{m-1}\mathcal B_\delta^jE_\delta 1
 &\le 2c\delta\bigl(w-\mathcal B_\delta^m w\bigr)
 \le 2c\delta w,\\
 \sum_{j=0}^{m-1}\mathcal B^jE_\delta 1
 &\le c\delta\bigl(\mathcal B^m1-1\bigr)
 \le c\delta(H-1).
 \end{aligned}
\]

Write $p=d/2$. For each $N\ge1$, set
\[
 L_N=\max_{1\le k\le N} k^p
          \sup_{x,y\in\R^d}\mathcal B_\delta^k(x,y).
\]
This quantity is finite since
$\mathcal B_\delta^k(x,y)\le e^{(b+\delta)k}\rho_k(x-y)$.
For $2\le n\le N$, expand
\[
 \mathcal B_\delta^n-\mathcal B^n
 =\sum_{j=0}^{n-1}\mathcal B_\delta^jE_\delta\mathcal B^{n-1-j}.
\]
For $j\le(n-1)/2$, \eqref{uu1} bounds the kernel of
$\mathcal B^{n-1-j}$ by $Cn^{-p}$. The contribution of these terms
at $(x,y)$ is at most
$Cn^{-p}\sum_{j=0}^{n-1}(\mathcal B_\delta^j E_\delta1)(x)$.
For $j>(n-1)/2$, the kernel of $\mathcal B_\delta^j$ is at most
$2^pL_N n^{-p}$. By symmetry of $E_\delta$ and $\mathcal B$,
the contribution of these terms is at most
$2^pL_N n^{-p}\sum_{j=0}^{n-1}(\mathcal B^jE_\delta1)(y)$.
The two telescoping estimates therefore give
\[
 0\le\mathcal B_\delta^n(x,y)-\mathcal B^n(x,y)
 \le C\delta(1+L_N)n^{-p}.
\]
For $n=1$, the same bound follows from
$E_\delta\le c\delta A\le c\delta e^b P$.
Using \eqref{uu1} and taking the maximum over $1\le n\le N$ gives
$L_N\le C_0+C\delta(1+L_N)$, with constants independent of
$N$ and $\delta$. Decreasing $\delta_0$ so that $C\delta_0\le1/2$
gives a uniform bound on $L_N$. Substitution into the preceding
display proves \eqref{eq:gap-kernel-stability}.
\end{proof}

\subsection{Renewal comparison and weak disorder beyond the second-moment threshold}

Retain $b=b_*$ and the rank-one kernel $\varepsilon f\otimes f$ from
Lemma~\ref{lem:rank-one-gauge}. For $n\ge1$, define
\[
 K(n):=\varepsilon\langle f,\mathcal B^{n-1}f\rangle,
 \qquad
 K_\delta(n):=\varepsilon\langle f,\mathcal B_\delta^{n-1}f\rangle .
\]
Write $M=\varepsilon f\otimes f$. In the splitting below, of the type introduced by Nummelin~\cite{Nummelin}, a unit path segment that stays in $C$ may receive a mark. The marked and unmarked segments contribute the kernels $M$ and $\mathcal B_\delta$, respectively. Two successive marks $n$ steps apart enclose $n-1$ unmarked steps. Integrating over the intermediate positions gives the rank-one identity
\[
 M\mathcal B_\delta^{n-1}M=K_\delta(n)M.
\]
The next lemma shows that $K$ is a probability law for the renewal gaps. For $\delta>0$ we use $K_\delta(n)$ as a weight and compare it with $K(n)$ uniformly in $n$.

\begin{lemma}\label{lem:marked-renewal}
The sequence $K$ satisfies $K(n)>0$ for all $n\ge1$,
$\sum_{n\ge1}K(n)=1$, and
\[
 K(n)\sim C_K n^{-d/2},\qquad
 C_K=\varepsilon(4\pi)^{-d/2}
       \left(\int_{\R^d} f(x)H(x)\,dx\right)^2>0.
\]
There exist $\delta_0>0$ and $C_0<\infty$ such that
\begin{equation}\label{eq:gap-weight-bound}
 K_\delta(n)\le (1+C_0\delta)K(n),
 \qquad n\ge1,\quad 0\le\delta\le\delta_0 .
\end{equation}
\end{lemma}
\begin{proof}
Write $J_f=\langle f,\psi\rangle>0$, where $\psi$ is the function
in \eqref{oo}. Since $S\psi=\psi$ by \eqref{ll},
$(I-\mathcal B)\psi=\varepsilon J_f f$. Multiplying by
$\mathcal B^j$, pairing with $f$, and summing gives
\[
 J_f\sum_{n=1}^N K(n)
   =\langle f,\psi\rangle-\langle f,\mathcal B^N\psi\rangle.
\]
For each $L>0$, \eqref{uu1} and $\mathcal B^N1\le H$ give
\[
 0\le(\mathcal B^N\psi)(x)
 \le C N^{-d/2}\int_{|y|\le L}\psi(y)\,dy
       +\|H\|_\infty\sup_{|y|>L}\psi(y).
\]
Let $N\to\infty$ and then $L\to\infty$. Since $\psi\in C_0(\R^d)$,
$\|\mathcal B^N\psi\|_\infty\to0$, and hence $\sum_{n\ge1}K(n)=1$.

At $n=1$, $K(1)=\varepsilon\|f\|_2^2>0$. For $n\ge2$,
$\mathcal B^{n-1}\ge P^{n-1}$ and the heat kernel is strictly positive.
The local limit \eqref{eq:regen-local-limit} gives the stated asymptotic
by dominated convergence against $f(x)f(y)\,dx\,dy$, using the
upper bound \eqref{uu1}. In particular, $K(n)\ge c n^{-d/2}$ for all
$n\ge1$. By Lemma~\ref{lem:gap-perturb}, for $n\ge2$,
\[
 0\le K_\delta(n)-K(n)
 \le C\varepsilon\delta\|f\|_1^2(n-1)^{-d/2}
 \le C_0\delta K(n).
\]
Since $K_\delta(1)=K(1)$, this proves \eqref{eq:gap-weight-bound}.
\end{proof}

Let $\widehat\tau$ be a renewal process starting at $0$ with gap
distribution $K$. Independently, let $(\xi_i)_{i\ge1}$ be independent
standard Gaussian variables. For $\sigma>0$ and $h\in\R$, set
\begin{equation}\label{eri}
r_i=e^{\sigma\xi_i-\sigma^2/2}
\end{equation}
and define the constrained partition functions
\[
 \mathcal U_0(h)=1,\qquad
 \mathcal U_n(h)=E_{\widehat\tau}\!\left[
 \mathbf1_{\{n\in\widehat\tau\}}
 \prod_{i\in\widehat\tau\cap[1,n]} e^h r_i\right],\qquad n\ge1.
\]
Here $E_{\widehat\tau}$ averages over the renewal process alone.

\begin{lemma}\label{lem:pinning-sum}
For every $\sigma>0$, there exists $h>0$ such that
$\sum_{n\ge0}\mathcal U_n(h)<\infty$ almost surely.
\end{lemma}
\begin{proof}
Let $a=d/2-1\ge1/2$ and $\varphi(n)=2\pi n^{1+a}K(n)$.
By Lemma~\ref{lem:marked-renewal}, $\varphi(n)\to2\pi C_K>0$.
Thus $\varphi$ is slowly varying and
\[
 \sum_{n\ge1}\frac{1}{n^{2(1-a)}\varphi(n)^2}=\infty.
\]
For $d=3$ the summand is asymptotic to a positive constant times $1/n$.
Berger and Lacoin~\cite[Theorem~2.2]{BL} therefore give
$h_c(\sigma)>0$ for the quenched critical reward in the normalization
$e^hr_i$. Fix $0<h<h_c(\sigma)$. For unnormalized contact weights
the reward is $h-\sigma^2/2$, strictly below the corresponding critical
point $h_c(\sigma)-\sigma^2/2$. Including a contact at time zero, as in
\cite{AZ}, multiplies each $\mathcal U_n(h)$ by the same strictly
positive random factor $e^{h+\sigma\xi_0-\sigma^2/2}$, where $\xi_0$
is an additional independent standard Gaussian. Mourrat's summability
theorem~\cite{Mourrat}, in the form of
\cite[Theorem~2.1(i)]{AZ}, now yields the assertion.
\end{proof}

\begin{theorem}[Strict separation of the thresholds]\label{thm:L2-ui}
There is a $\delta_*>0$ such that $(Z_t^\kappa(0))_{t\ge0}$ is
uniformly integrable for every $\kappa\in[b_*,b_*+\delta_*]$.
In particular, $b_*<\kappa_c$.
\end{theorem}
\begin{proof}
Keep the ball $C\supset\operatorname{supp} f$ from
Lemma~\ref{lem:rank-one-gauge} and fix
\[
 \sigma=\sqrt b\inf_{x\in C}q(x)>0 .
\]
Choose $h>0$ as in Lemma~\ref{lem:pinning-sum}, and then
$0<\delta_*\le\delta_0$ such that $\log(1+C_0\delta_*)\le h$ in
\eqref{eq:gap-weight-bound}. Fix $\kappa=b+\delta$ with
$0\le\delta\le\delta_*$. The bridge formula \eqref{ebridge} for $S_\kappa^C$, Lemma~\ref{lem:rank-one-gauge}\textup{(4)}, and $q\ge0$ give
$\varepsilon f(x)f(y)\le S_b^C(x,y)\le S_\kappa^C(x,y)$.

For each $n\ge1$, equip paths $D$ starting at zero and auxiliary
coordinates $u=(u_1,\ldots,u_n)\in[0,1]^n$ with the finite measure
\[
 \mu_n(dD,du)=e^{\kappa\int_0^n q(D_s)\,ds}P_0^D(dD)\,du .
\]
Here $du$ is product Lebesgue measure. Write $m_i(D,u)=1$ when
$D_s\in C$ for all $s\in[i-1,i]$ and
\[
 u_i<\frac{\varepsilon f(D_{i-1})f(D_i)}{S_\kappa^C(D_{i-1},D_i)} .
\]
Set $m_i=0$ otherwise and interpret the ratio as zero if its denominator
vanishes. Integrating the Feynman--Kac weight on one marked segment
gives the kernel $\varepsilon f\otimes f$, since conditional integration
over $u_i$ multiplies $S_\kappa^C(x,y)$ by the displayed ratio. The
unmarked segment therefore has kernel $\mathcal B_\delta$.

With $(r_i)_{i\ge1}$ as in \eqref{eri}, we let
\[
 G_n=\int\prod_{i=1}^n r_i^{m_i(D,u)}\,\mu_n(dD,du).
\]
The contribution from paths with no marks is $(\mathcal B_\delta^n1)(0)$.
For paths whose marked times are precisely $1\le t_1<\cdots<t_k\le n$,
the Markov property and the rank-one formula give the contribution
\begin{equation}\label{eq:gap-marked-product}
 \varepsilon(\mathcal B_\delta^{t_1-1}f)(0)
 \left(\prod_{j=2}^k K_\delta(t_j-t_{j-1})\right)
 \langle f,\mathcal B_\delta^{n-t_k}1\rangle
 \prod_{j=1}^k r_{t_j}.
\end{equation}
The powers of order zero act as the identity. By
Lemma~\ref{lem:gap-perturb}, the no-mark contribution and the terminal
factor $\langle f,\mathcal B_\delta^{n-t_k}1\rangle$ are bounded
uniformly. The same lemma and \eqref{uu1} bound the first factor by
$C_1(t_1-1)^{-d/2}$ for $t_1\ge2$, hence by $C_2K(t_1)$ by
Lemma~\ref{lem:marked-renewal}. At $t_1=1$ this last bound follows
from $K(1)>0$. For each subsequent gap,
\eqref{eq:gap-weight-bound} gives $K_\delta(t_j-t_{j-1})\le
e^h K(t_j-t_{j-1})$. Grouping \eqref{eq:gap-marked-product} by
the last mark $t_k$ and using the definition of $\mathcal U_{t_k}(h)$,
we obtain
\begin{equation}\label{eq:critical-random-bound}
 G_n\le C_3\sum_{j=0}^n\mathcal U_j(h)
 \le C_3\sum_{j=0}^\infty\mathcal U_j(h)=:\mathcal R_h<\infty
 \quad\text{almost surely},
\end{equation}
with $C_3$ independent of $n$ and $\delta$. Only $K$ serves as the
renewal law in this bound.

Let $\eta$ be a standard Brownian motion independent of $D$. For fixed
$z=(D,u)$ define the centered Gaussian fields
\[
 X(z)=\sqrt\kappa\int_0^n q(D_s)\,d\eta_s,
 \qquad Y(z)=\sigma\sum_{i=1}^n\xi_i m_i(z).
\]
When two paths $z=(D,u)$ and $z'=(D',u')$ share a mark at $i$, both
remain in $C$ throughout $[i-1,i]$. Since $q\ge0$ and $\kappa\ge b$,
\begin{equation}\label{eq:critical-covariance-order}
 \operatorname{Cov}(X(z),X(z'))
 =\kappa\int_0^n q(D_s)q(D_s')\,ds
 \ge\sigma^2\sum_{i=1}^n m_i(z)m_i(z')
 =\operatorname{Cov}(Y(z),Y(z')).
\end{equation}
The integral of $e^{X-\operatorname{Var}(X)/2}$ against $\mu_n$ is
$F_n^\kappa(\eta)$ from Proposition~\ref{prop:scalar}. The corresponding
integral for $Y$ is $G_n$, since $m_i^2=m_i$. For $\lambda>0$, let
$\Phi_\lambda(x)=1-e^{-\lambda x}$. By \eqref{eq:critical-covariance-order}
and the concavity of $\Phi_\lambda$, Kahane's inequality
\cite[Lemma~1]{Kahane} gives
\[
 \E\Phi_\lambda(F_n^\kappa)\le\E\Phi_\lambda(G_n).
\]
Indeed, sample independent paths from $\mu_n/\mu_n(1)$, independently
of the Gaussian fields, and apply the finite-sum inequality to the
empirical integrals with equal weights summing to $\mu_n(1)$. Both
integrals have finite expectation $\mu_n(1)\le e^{\kappa n}$, so the
conditional strong law gives convergence of the empirical integrals
almost surely. Bounded convergence then proves the inequality.

Taking $\lambda=1/L$ and using \eqref{eq:critical-random-bound}, we obtain
\[
 \sup_{n\ge1}\Pp(F_n^\kappa>L)
 \le\frac{\E[1-e^{-\mathcal R_h/L}]}{1-e^{-1}}
 \longrightarrow0\qquad (L\to\infty).
\]
Proposition~\ref{prop:scalar} yields uniform integrability for every
$\kappa\in[b_*,b_*+\delta_*]$. By
\eqref{eq:critical-parameter} and Proposition~\ref{prop:mean-dichotomy},
$\kappa_c\ge b_*+\delta_*>b_*$.
\end{proof}

\begin{proof}[Proof of Theorem~\ref{main}]
Lemma~\ref{lem:spectral}, Theorem~\ref{thm:L2-ui}, and
Proposition~\ref{prop:critical-parameter} give
$0<b_*<\kappa_c<\infty$. Since $\kappa_c>b_*$,
Theorem~\ref{thm:gap-case} gives uniform integrability at $\kappa_c$.
By Propositions~\ref{prop:mean-dichotomy} and~\ref{prop:critical-parameter},
$c_\kappa=1$ for every $0\le\kappa\le\kappa_c$. The Laplace-functional
limit in \eqref{eq:pi-laplace} and
Proposition~\ref{prop:mean-dichotomy} now give
\[
 \mathcal L_m^\kappa(X_t)\Rightarrow\pi_\kappa,\qquad
 \int_{\Mloc}\mu(f)\,\pi_\kappa(d\mu)=\int_{\R^d}f(x)\,dx,
 \quad f\in C_c^+(\R^d).
\]
Taking a nonzero $f\in C_c^+(\R^d)$ shows that
$\pi_\kappa\ne\delta_0$. This proves \textup{(1)}.

For $\kappa>\kappa_c$,
Propositions~\ref{prop:critical-parameter} and~\ref{prop:mean-dichotomy}
give $c_\kappa=0$ and $Z_t^\kappa(0)\to0$ almost surely.
As $\kappa>\kappa_c>b_*$, Theorem~\ref{thm:sharp} gives, for every
$\theta\in(0,1)$,
\[
 \E[(Z_t^\kappa(0))^\theta]\le Ce^{-ct},\qquad t\ge0.
\]
Thus \eqref{fracdecay} holds, and Proposition~\ref{prop:eventual-extinction}
proves \textup{(2)}. The characterization of $\kappa_c$ follows from
\eqref{eq:critical-parameter} and the equivalence in
Proposition~\ref{prop:mean-dichotomy} between $c_\kappa=1$ and
uniform integrability of $(Z_t^\kappa(0))_{t\ge0}$.
\end{proof}

\begin{proof}[Proof of Theorem~\ref{thm:critical-intro}]
Theorem~\ref{thm:L2-ui} gives $b_*<\kappa_c$, and
Theorem~\ref{thm:gap-case} therefore applies at $\kappa_c$.
Together with Propositions~\ref{prop:critical-parameter}
and~\ref{prop:mean-dichotomy}, this gives uniform integrability and
$\E Z_\infty^\kappa=1$ throughout $[b_*,\kappa_c]$.
Proposition~\ref{prop:exact-L2} gives
$\lim_{t\to\infty}\E[(Z_t^\kappa)^2]=\infty$ for each such $\kappa$.

For $0\le\kappa\le\kappa_c$, positivity of $\E Z_\infty^\kappa$
implies $\Pp(Z_\infty^\kappa>0)>0$. By Lemma~\ref{lem:free-energy},
$p(\kappa)\le0$, and a strictly negative value would force
$Z_t^\kappa\to0$ almost surely. Thus $p(\kappa)=0$ on this interval.
For $\kappa>\kappa_c$, Propositions~\ref{prop:critical-parameter}
and~\ref{prop:mean-dichotomy} give strong disorder. As $\kappa>b_*$,
Theorem~\ref{thm:sharp} and Lemma~\ref{lem:free-energy}
give $p(\kappa)<0$. The almost sure free-energy limit and the
fractional-moment limits at $\kappa_c$ follow from
Proposition~\ref{prop:critical-free-energy}.\\
\end{proof}

\bigskip

\noindent  {\bf  \large  AI Declaration} \\
 
The authors benefited from discussions with commercial large language models. The authors obtained the results, worked out the principal proofs, and take full responsibility for the final manuscript.

\bibliographystyle{amsplain}

\end{document}